# Homologie et cohomologie de Hochschild de certaines algèbres polynomiales classiques et quantiques

Lionel Richard
*Laboratoire de Mathématiques Pures*
*Université Blaise Pascal, 63177 Aubière Cedex, France*
E-mail : Lionel.Richard@math.univ-bpclermont.fr

(Version provisoire : 25/06/2002)

**Résumé.**
Nous étudions l'homologie et la cohomologie de Hochschild de certaines classes d'algèbres polynomiales de type à la fois classique et quantique. Nous démontrons en particulier que l'algèbre des opérateurs différentiels sur tout espace quantique multiparamétré admet l'homologie de Hochschild d'une algèbre de Weyl usuelle, et satisfait une dualité de Poincaré.

## Introduction.

Le but de cet article est d'obtenir des calculs explicites sur l'homologie et la cohomologie de Hochschild d'une classe d'algèbres non-commutatives de polynômes, à la fois quantiques (au sens où elles contiennent comme sous-algèbres des copies du plan quantique de Manin) et classiques (au sens où elles contiennent aussi des copies de l'algèbre de Weyl $A_1$). Plus précisement, sur un corps de base $k$ commutatif de caractéristique nulle, les algèbres considérées sont paramétrées par un entier $n \geq 1$, une matrice multiplicativement antisymétrique $\Lambda = (\lambda_{i,j})_{1\leq i,j\leq n}$ à coefficients non-nuls dans $k$, et un entier $0 \leq r \leq n$. On note $A_{n,r}^{\Lambda}$ l'algèbre engendrée sur $k$ par $n+r$ générateurs $y_1, \ldots, y_n, x_1, \ldots, x_r$ tels que :

(1) la sous-algèbre de $A_{n,r}^{\Lambda}$ engendrée par $y_1, \ldots, y_n$ est isomorphe à l'espace quantique $\mathcal{O}_\Lambda(k^n)$, c'est-à-dire que $y_iy_j = \lambda_{i,j}y_jy_i$ pour tous $1 \leq i,j \leq n$ ;

(2) si $r \geq 1$, la sous-algèbre de $A_{n,r}^{\Lambda}$ engendrée par $x_i$ et $y_i$ est pour tout $1 \leq i \leq r$ isomorphe à l'algèbre de Weyl $A_1(k)$, c'est-à-dire que $x_iy_i - y_ix_i = 1$ ;

(3) les autres relations entre générateurs sont de la forme $x_ix_j = \lambda_{i,j}x_jx_i$ pour tous $1 \leq i,j \leq r$, et $x_iy_j = \lambda_{i,j}^{-1}y_jx_i$ pour tous $1 \leq i \leq r,\ 1 \leq j \leq n,\ i \neq j$.

Ce type d'algèbres, qui s'interprètent naturellement comme opérateurs différentiels tordus sur l'espace quantique $\mathcal{O}_\Lambda(k^n)$, intervient dans de nombreuses situations en liaison avec l'étude des algèbres de Weyl quantiques, des extensions de Ore itérées, et de l'équivalence birationnelle de telles algèbres. On donnera dans la première section de l'article des précisions sur ce point. Contentons-nous de remarquer ici que, dans le cas purement quantique où $r = 0$, l'algèbre $A_{n,r}^\Lambda$ se réduit à l'espace quantique $\mathcal{O}_\Lambda(k^n)$, et que dans le cas purement classique où $r = n$ et tous les $\lambda_{i,j}$ valent 1, l'algèbre $A_{n,n}^{(1)}$ n'est autre que l'algèbre de Weyl usuelle $A_n(k) = A_1(k) \otimes \cdots \otimes A_1(k)$. L'homologie de Hochschild dans ces deux situations est bien connue (dans le cas classique, voir [26], [17] ; pour le cas quantique, [27], [14] ou [12]). Les deux situations générales qui font l'objet principal de cet article concernent quant à elles d'une part le cas dit "libre" où les coefficients de la matrice $\Lambda$ forment dans $k^*$ un groupe abélien libre de rang maximal, et d'autre part le cas dit "semi-classique" des opérateurs différentiels usuels sur l'espace quantique $\mathcal{O}_\Lambda(k^n)$ (c'est-à-dire le cas où $r = n$). Nous montrons dans ce dernier cas que l'homologie et la cohomologie de Hochschild de $A_{n,n}^\Lambda$ sont, sans aucune hypothèse restrictive sur la matrice $\Lambda$ de quantification de l'espace, celles de l'algèbre de Weyl usuelle sur un espace de fonctions polynomiales commutatif.

La première section de l'article est consacrée d'une part aux préliminaires nécessaires concernant les algèbres $A_{n,r}^\Lambda$, et d'autre part à une présentation de ces algèbres comme cas particuliers d'algèbres enveloppantes quantiques définies par M. Wambst. D'une résolution libre introduite pour de telles algèbres en [27], on déduit un complexe de Koszul quantique dont l'homologie est celle des algèbres $A_{n,r}^\Lambda$. Nous montrons à la section 2 que l'on peut extraire de ce complexe de Koszul $(K, d)$ un sous-complexe $(K_C, d)$ plus petit, tel que le complexe quotient $K/K_C$ soit acyclique. Ce petit complexe nous permet dans la suite le calcul explicite de l'homologie de Hochschild dans deux situations particulières significatives : le cas "semi-classique" et le cas "libre". Le cas "semi-classique" fait l'objet de la section 3 : on y démontre que, pour tout $n \geq 1$ et quelle que soit la matrice $\Lambda$, l'homologie de Hochschild de $A_{n,n}^\Lambda$ est concentrée en degré $2n$, où elle vaut $k$. A la section 4, sous l'hypothèse que les coefficients $\Lambda$ forment un groupe abélien libre de rang maximal $n(n-1)/2$, on calcule pour tout $1 \leq r < n$ l'homologie de Hochschild de $A_{n,r}^\Lambda$, qui est cette fois concentrée en degrés 0, 1 et $2r$. Le plus petit cas (en termes de dimension) qui échappe aux hypothèses précédemment traitées est celui où $n = 2$ et $r = 1$, avec l'unique coefficient $\lambda$ de la matrice $\Lambda$ qui est racine de l'unité dans $k$ ; l'homologie de Hochschild de l'algèbre $A_{2,1}^\Lambda$ correspondante est calculée à la section 5. La dualité de Poincaré pour les algèbres semi-classiques $A_{n,n}^\Lambda$ est le théorème central de la section 7. Sa preuve repose sur les résultats obtenus à la section 6 sur la cohomologie de Hochschild des algèbres $A_{n,r}^\Lambda$ en général. Enfin, les remarques et résultats exploratoires regroupés dans la section 8 montrent, au vu de ce qui a été établi à la section 5, que l'on ne peut pas espérer dans ce contexte général de propriété aussi favorable de dualité.

# 1 Préliminaires.

Dans tout le texte $k$ désigne un corps commutatif de caractéristique nulle. On note $\otimes$ le produit tensoriel sur $k$ pour des espaces vectoriels : $\otimes = \otimes_k$. Si $V$ et $W$ sont deux $k$-espaces vectoriels, et $f$ une application linéaire de $V \otimes V$ dans $W$, on utilise la convention suivante : pour tout $n \geq 2$ et pour tout $0 \leq i \leq n-1$ on note

$$f_i : V^{\otimes n} \to V^{\otimes (i-1)} \otimes W \otimes V^{\otimes n-i-1}$$

l'application linéaire définie par : $f_i = \mathrm{id}_{V^{\otimes i-1}} \otimes f \otimes \mathrm{id}_{V^{\otimes n-i-1}}$. Si $n$ est un entier non nul, pour tout $i \leq n$ on note $[i] = (0, \dots, 0, 1, 0, \dots, 0)$ le $i^{\text{ème}}$ vecteur de la base canonique de $\mathbb{N}^n$. Pour un multi-entier $\gamma \in \mathbb{N}^n$, on note $|\gamma| = \gamma_1 + \dots + \gamma_n$ la longueur de $\gamma$. Enfin pour $\alpha \in \mathbb{N}^r$ et $\beta \in \mathbb{N}^n$, on note globalement le monôme $x^\alpha y^\beta = x_1^{\alpha_1} \dots x_r^{\alpha_r} y_1^{\beta_1} \dots y_n^{\beta_n}$.
L'homologie et la cohomologie d'une algèbre sont dans ce texte toujours calculées à valeurs dans l'algèbre elle-même.

## 1.1 Définition des algèbres $A_{n,r}^\Lambda$.

Soient $r$ et $n$ deux entiers tels que $n \geq 1$ et $0 \leq r \leq n$, et soit $\Lambda = (\lambda_{i,j}) \in M_n(k)$ une matrice multiplicativement antisymétrique de taille $n$, c'est-à-dire telle que pour tous $i, j$ on ait : $\lambda_{i,j}\lambda_{j,i} = \lambda_{i,i} = 1$. On note $A_{n,r}^\Lambda$ la $k$-algèbre engendrée par $n+r$ générateurs $x_1, \dots, x_r, y_1, \dots, y_r, y_{r+1}, \dots, y_n$ soumis aux relations suivantes :

$$\begin{array}{lll} y_i y_j = \lambda_{i,j} y_j y_i & \text{si} & 1 \leq i, j \leq n, \\ x_i y_j = \lambda_{i,j}^{-1} y_j x_i & \text{si} & 1 \leq i \leq r,\ 1 \leq j \leq n,\ i \neq j, \\ x_i y_i = y_i x_i + 1 & \text{si} & 1 \leq i \leq r, \\ x_i x_j = \lambda_{i,j} x_j x_i & \text{si} & 1 \leq i, j \leq r. \end{array}$$

## 1.2 Remarques.

**1.2.1** Il est facile de voir que l'algèbre $A_{n,r}^\Lambda$ est une extension itérée de Ore. Elle est noethérienne et intègre, de dimension de Gelfand-Kirillov égale à $n + r$, et les monômes $\{x_1^{\alpha_1} \dots x_r^{\alpha_r} y_1^{\beta_1} \dots y_n^{\beta_n}\}_{\alpha \in \mathbb{N}^r, \beta \in \mathbb{N}^n}$ en forment une base de $k$-espace vectoriel. (Pour des propriétés générales des algèbres $A_{n,r}^\Lambda$, voir [25]).

**1.2.2** On peut montrer que toute algèbre associative engendrée par $N$ générateurs et des relations du type $xy = \lambda yx$ ou $xy = yx + 1$, et possédant une base de PBW, est isomorphe à une algèbre $A_{n,r}^\Lambda$ avec $n + r = N$ (voir [25]).

**1.2.3** Les algèbres $A_{n,r}^\Lambda$ peuvent être vues comme des algèbres d'opérateurs différentiels sur un espace quantique paramétré par $\Lambda$. Notons $\mathcal{O}_\Lambda(k^n)$ l'algèbre engendrée par les $y_i$ ; c'est un espace affine quantique. Pour tout $i$ on note $\mu_i$ la multiplication à *gauche* par $y_i$ dans $\mathcal{O}_\Lambda(k^n)$, et pour $j \leq r$ on définit $\partial_j \in \mathrm{End}(\mathcal{O}_\Lambda(k^n))$ par :

$$\partial_j : y_1^{m_1} \dots y_n^{m_n} \mapsto m_j \prod_{k<j} (\lambda_{k,j})^{m_k} y_1^{m_1} \dots y_j^{m_j - 1} \dots y_n^{m_n}.$$

C'est une $\sigma_j$-dérivation, où $\sigma_1, \ldots, \sigma_r$ sont les automorphismes de $\mathcal{O}_\Lambda(k^n)$ définis par :

$$\sigma_j : P(y_1, \ldots, y_n) \mapsto P(\lambda_{1,j}\ y_1, \ldots, \lambda_{n,j}\ y_n).$$

On vérifie alors que la sous-algèbre de $\mathrm{End}(\mathcal{O}_\Lambda(k^n))$ engendrée par les $\partial_j$ et les $\mu_i$ est isomorphe à $A_{n,r}^\Lambda$.

## 1.3 Exemples et terminologie.

**1.3.1** *Cas purement quantique.* Dans le cas où $r = 0$, l'algèbre $A_{n,0}^\Lambda$ n'est autre que l'espace affine quantique $\mathcal{O}_\Lambda(k^n)$, dont l'homologie de Hochschild est calculée au théorème 6.1 de [27], ainsi que dans [14] et [12].

**1.3.2** *Cas purement classique.* Dans le cas où $r = n$ et où tous les $\lambda_{i,j}$ valent 1, l'algèbre $A_{n,n}^{(1)}$ est l'algèbre de Weyl usuelle $A_n(k)$, dont l'homologie est bien connue (voir par exemple [26] ou [17]), et rappelée au paragraphe 3.1.

**1.3.3** *Cas semi-classique.* Dans le cas $n = r$, sans hypothèse sur la matrice $\Lambda$, les algèbres $A_{n,n}^\Lambda$ sont des cas particuliers d'algèbres de Weyl quantiques $A_n^{\bar{q},\Lambda}(k)$ (voir plus loin au 1.4.2), correspondant au cas où les paramètres de quantification $\bar{q} = (q_1, \ldots, q_n)$ valent tous 1. Elles apparaissent par exemple en [11] (Exemple 2.1, lorsque $q^2 = 1$) et dans [9]. Il s'agit également des algèbres d'opérateurs différentiels sur un espace quantique présentées par M. Wambst dans la section 11 de [27], lorsque $q = 1$ et $p_i = 1$ pour tout $i$, mais sans que soit fait le calcul explicite de leur homologie. Enfin, lorsque tous les $\lambda_{i,j}$ valent $-1$ pour $i \neq j$, ces algèbres apparaissent en physique mathématique comme algèbres de “pseudo-bosons” (voir [5]). Les algèbres $A_{n,n}^\Lambda$ sont simples, et de centre $k$ (voir [11], [24]). Pour ces algèbres, nous avons calculé dans [24] les modules de degré 0 et 1 en cohomologie de Hochschild, et de degré 0 en homologie, et constaté qu'ils étaient égaux à ceux de $A_n(k)$. On verra ici que toute l'homologie et la cohomologie de Hochschild des $A_{n,n}^\Lambda$ sont identiques à celles de l'algèbre de Weyl $A_n(k)$.

**1.3.4** *Cas libre.* Sans faire d'hypothèse sur les entiers $r$ et $n$, on s'intéresse au cas où les $\lambda_{i,j}$ forment un groupe abélien libre de rang maximal $n(n-1)/2$. Quoiqu'assez restrictive, cette hypothèse se retrouve fréquemment dans la littérature concernant les groupes quantiques (c'est par exemple le cas des “general quantum functions” de V. Artamonov dans [3]).

**1.3.5** Si tous les $\lambda_{i,j}$ valent 1, l'algèbre $A_{n,r}^{(1)}$ est le produit tensoriel de l'algèbre de Weyl classique $A_r(k)$ et de l'algèbre des polynômes commutatifs en $n-r$ variables. Le calcul de l'homologie de Hochschild se ramène alors au cas 1.3.2 grâce à la formule de Künneth.

**1.3.6** *Cas mixte minimal.* Le premier cas échappant aux cas généraux précédents s'obtient en prenant $r = 1$, $n = 2$ et un coefficient $\lambda$ d'ordre fini supérieur ou égal à 2 dans $k^*$. Sans hypothèse sur $\lambda$, le cas $n = 2, r = 1$ sera appelé “mixte minimal”. Des calculs explicites d'homologie et de cohomologie seront faits pour ce cas aux sections 5 et 8.

## 1.4 Autres exemples.

**1.4.1** Les algèbres $A_{n,r}^{\Lambda}$ sont un cas particulier des algèbres $S_{\mathbf{Q},\mathbf{P}}^{p}$ définies par J.A. et J.J. Guccionne dans [13], où $p = n$, c'est-à-dire que l'on ne localise par rapport à aucune variable. Cependant pour $r \geq 1$ les auteurs ne font pas le calcul explicite de l'homologie dans le cadre que nous considérons ici. Les algèbres $A_{n,r}^{\Lambda}$ peuvent également être vues comme des cas particuliers des algèbres $Q$-résolubles définies par A. Panov dans [22].

**1.4.2** Beaucoup d'exemples d'algèbres liées aux groupes quantiques ont des localisations communes avec les algèbres $A_{n,r}^{\Lambda}$. D.A. Jordan étudie dans [16] une algèbre $S_{p,q}$ d'opérateurs différentiels eulériens sur l'algèbre de polynômes $k[y]$, qui lorsque le paramètre $p$ vaut 1, admet une localisation commune avec une algèbre $A_{2,1}^{\Lambda}$. La version alternative de l'algèbre de Heisenberg quantique présentée à la section 2 de [18] admet également une localisation commune avec $A_{2,1}^{\Lambda}$. Par ailleurs, les algèbres de Weyl quantiques multiparamétrées $A_n^{\bar{q},\Lambda}(k)$ (voir la définition I.2.6 de [6]) introduites par G.Maltsiniotis dans [20], et qui ont fait l'objet de nombreuses études du point de vue de la théorie des anneaux (voir par exemple [1], [15], [11], [9], [10] et leurs bibliographies) admettent toujours une localisation commune avec des algèbres $A_{n,r}^{\Lambda}$, en particulier si certains des paramètres $q_i$ valent 1 (voir [24]). Dans le cas purement quantique, on peut sur ce point se référer aux travaux de D.A. Jordan ([15]), de G. Cauchon ([7]) et d'A. Panov ([22]).

**1.4.3** Plus généralement les algèbres $A_{n,r}^{\Lambda}$ jouent un rôle central dans un certain nombre de questions d'équivalence rationnelle non-commutative (voir [25]). En particulier pour $r = n = 2$ le corps de fractions de $A_{2,2}^{\Lambda}$ est étudié dans [2] comme contre-exemple à une version mixte de la conjecture de Gelfand-Kirillov. Enfin le théorème 4.2 de [23] traite du lien entre l'équivalence rationnelle des algèbres $A_{n,n}^{\Lambda}$ et les tores quantiques multiparamétrés.

**1.4.4** Les $A_{n,r}^{\Lambda}$ sont des exemples particuliers de la notion d'algèbres enveloppantes généralisées d'algèbres de Lie quantiques introduites par M. Wambst dans [27]. C'est dans ce cadre que nous nous plaçons au paragraphe suivant, afin de décrire un complexe permettant un calcul explicite de l'homologie de Hochschild de $A_{n,r}^{\Lambda}$.

## 1.5 Homologie de Hochschild de certaines algèbres enveloppantes généralisées d'algèbres de Lie quantiques.

Dans ce paragraphe nous explicitons le cadre homologique adapté aux algèbres $A_{n,r}^{\Lambda}$.

**Définition 1.5.1** *Soit $V$ un $k$-espace vectoriel de dimension finie, et de base $(v_1, \ldots, v_m)$. Soient $Q = (q_{i,j}) \in M_m(k^*)$ matrice multiplicativement antisymétrique, et $f$ une forme linéaire définie sur $V \otimes V$ telle que*

$$f(v_i \otimes v_j) = -q_{i,j} f(v_j \otimes v_i).$$

*On appelle algèbre enveloppante quantique associée à ces données l'algèbre notée $U_Q(V,f)$, quotient de l'algèbre tensorielle sur $V$ par l'idéal engendré par la famille $\{v_i \otimes v_j - q_{i,j} v_j \otimes v_i - f(v_i, v_j),\ 1 \le i < j \le m\}$.*

**Remarques :** • L'algèbre enveloppante quantique $U_Q(V,f)$ est un cas particulier d'algèbre enveloppante généralisée $\mathcal{U}_c(\mathfrak{g}, f)$ d'algèbre de Lie quantique introduite par M.Wambst ([27], déf. 9.1). Rappelons pour mémoire qu'une symétrie de Hecke de marque $\nu \in k^*$ sur un $k$-espace vectoriel $V$ est une application linéaire $c : V \otimes V \to V \otimes V$ telle que :

$$c_1c_2c_1 = c_2c_1c_2, \text{ et } (c-\nu)(c+1) = 0.$$

Nous nous restreignons ici au cas où l'algèbre de Lie est abélienne et où la symétrie de Hecke $c : V \otimes V \to V \otimes V$ se réduit à des $q$-commutations régies par la matrice $Q$. Ceci correspond à l'exemple 1.2(5) de [27], dans le cas où la marque $q$ de $c$ vaut 1.
Nous noterons dans la suite $c(Q)$ cette symétrie de Hecke particulière.

• Pour $f = 0$ on retrouve la notion d'algèbre symétrique quantique, telle qu'elle est définie dans [27] en 2.1. En effet $U_Q(V,0) = S_QV$ est le quotient de l'algèbre tensorielle $T(V)$ par l'idéal engendré par les éléments $v_i \otimes v_j - q_{i,j} v_j \otimes v_i$. C'est exactement l'espace affine quantique $\mathcal{O}_Q(k^m)$.

• De manière analogue on définit l'algèbre extérieure quantique ([27], déf. 2.1) $\Lambda_Q V$ comme le quotient de l'algèbre tensorielle $T(V)$ par l'idéal engendré par la famille $\{v_i \otimes v_j + q_{i,j} v_j \otimes v_i,\ i < j\}$.

M. Wambst donne dans l'article [27], sous certaines conditions, un complexe de Koszul quantique dont l'homologie est l'homologie de Hochschild de l'algèbre $\mathcal{U}_c(\mathfrak{g}, f)$. Explicitons-le pour les algèbres $U_Q(V,f)$. Fixons un $k$-espace vectoriel $V$ de dimension finie, une base $(v_1, \dots, v_m)$ de $V$, une matrice $Q \in M_m(k^*)$ multiplicativement antisymétrique, et une application linéaire $f$ de $V \otimes V$ dans $k$ vérifiant la condition de la définition 1.5.1. L'espace vectoriel $\Lambda_Q V$ est un espace $\mathbb{N}$-gradué, et il a pour base la famille $(v_1^{\gamma_1} \wedge \dots \wedge v_m^{\gamma_m})_{\gamma \in \{0,1\}^m}$. Le degré d'un monôme $v_1^{\gamma_1} \wedge \dots \wedge v_m^{\gamma_m}$ est égal à $|\gamma|$. Pour tout $* \in \mathbb{N}$, on note $\Lambda_Q^* V$ la partie homogène de degré $*$ de $\Lambda_Q V$. On gradue alors l'espace $U_Q(V,f) \otimes \Lambda_Q V \otimes U_Q(V,f)$ par le degré de $\Lambda_Q V$, et on en fait un complexe différentiel, à l'aide de la différentielle $\partial$ définie pour $a, b \in U_Q(V,f)$ par :

$$\begin{aligned}\partial(a \otimes v_{i_1} \wedge \dots \wedge v_{i_*} \otimes b) &= \\ \sum_{k=1}^{*} (-1)^{k-1} \Bigg( \Big( \prod_{s<k} q_{i_s, i_k} \Big) a v_{i_k} \otimes v_{i_1} \wedge \dots \widehat{v}_{i_k} \dots \wedge v_{i_*} \otimes b \\ - \Big( \prod_{s>k} q_{i_k, i_s} \Big) a \otimes v_{i_1} \wedge \dots \widehat{v}_{i_k} \dots \wedge v_{i_*} \otimes v_{i_k} b \Bigg)\end{aligned} \tag{1}$$

où $v_{i_1} \wedge \dots \widehat{v}_{i_k} \dots \wedge v_{i_*}$ désigne le produit extérieur $v_{i_1} \wedge \dots \wedge v_{i_*}$ auquel on a ôté le terme $v_{i_k}$.

**Proposition 1.5.2** *Soient $V$, $Q$ et $f$ comme dans la définition 1.5.1. Notons $U = U_Q(V, f)$, et $U^e = U \otimes U^{op}$. On suppose que le gradué associé à $U$ pour la filtration naturelle est l'algèbre symétrique $S_Q V$. Pour tout $2 \leq k \leq n$ posons $\Pi_k = c(Q)_1 \circ \ldots c(Q)_{k-1}$, et $\Pi_1 = Id$. De même, pour tout $1 \leq k \leq n-1$ posons $\check{\Pi}_k = c(Q)_{n-1} \circ \ldots c(Q)_k$, et $\check{\Pi}_n = Id$. Supposons en outre que pour tout $n \geq 2$ l'application $f'$ définie de $V^{\otimes n}$ dans $V^{\otimes n-2}$ par :*

$$f' = \sum_{1 \leq i < j \leq n} (-1)^{i+j+1} \left( (f \otimes I_{n-2})(I_1 \otimes \Pi_{j-1})\Pi_i - (I_{n-2} \otimes f)(\check{\Pi}_i \otimes I_1)\check{\Pi}_j \right),$$

*soit nulle. Alors le complexe $(U \otimes \Lambda_Q V \otimes U, \partial)$ est une résolution libre de $U$ par des $U^e$-bimodules.*

*Preuve.* C'est la proposition 10.3 de [27] dans le cas où $\mathcal{U}_c(\mathfrak{g}, f) = U_Q(V, f)$. $\square$

Ainsi l'homologie de Hochschild de $U = U_Q(V, f)$ est donnée par l'homologie du complexe $U \otimes_{U^e} (U \otimes \Lambda_Q^* V \otimes U)$ muni de la différentielle $\mathrm{id} \otimes \partial$. On identifie alors canoniquement $U \otimes_{U^e} (U \otimes \Lambda_Q^* V \otimes U)$ à $U \otimes \Lambda_Q^* V$, ce qui donne le complexe suivant, noté $(K(U)_*, d)$. Pour tous $\alpha \in \mathbb{N}^m$, $\gamma \in \{0, 1\}^m$, notons

$$v^\alpha \otimes v^\gamma = v^\alpha \otimes v_1^{\gamma_1} \wedge \ldots \wedge v_m^{\gamma_m}.$$

Comme $k$-espace vectoriel,

$$K(U) = \bigoplus_{\substack{\alpha \in \mathbb{N}^m \\ \gamma \in \{0,1\}^m}} k.v^\alpha \otimes v^\gamma,$$

le degré différentiel d'un monôme $v^\alpha \otimes v^\gamma$ étant égal à $|\gamma|$. La différentielle $d$ est donnée par :

$$d(v^\alpha \otimes v^\gamma) = \sum_{i=1}^m \Omega(\alpha, \gamma; i) v^\alpha v_i \otimes v^{\gamma - [i]} + \sum_{i=1}^m \Theta(\alpha, \gamma; i) v_i v^\alpha \otimes v^{\gamma - [i]}, \tag{2}$$

avec les coefficients :

$$\begin{aligned} &\Omega(\alpha, \gamma; i) = (-1)^{\sum_{k<i} \gamma_k} \prod_{k<i} q_{k,i}^{\gamma_k} \quad \text{si } \gamma_i = 1, \\ &\text{et} \\ &\Omega(\alpha, \gamma; i) = 0 \quad \text{si } \gamma_i = 0; \end{aligned} \tag{3}$$

pour le premier terme de la somme, et

$$\begin{aligned} &\Theta(\alpha, \gamma; i) = (-1)^{|\gamma| + \sum_{k>i} \gamma_k} \prod_{k>i} q_{i,k}^{\gamma_k} \quad \text{si } \gamma_i = 1, \\ &\text{et} \\ &\Theta(\alpha, \gamma; i) = 0 \quad \text{si } \gamma_i = 0. \end{aligned} \tag{4}$$

pour le deuxième.

**Proposition 1.5.3** *Reprenons les hypothèses de la proposition 1.5.2. Alors l'homologie de Hochschild de l'algèbre $U = U_Q(V, f)$ est l'homologie du complexe $(K(U)_*, d)$ que l'on vient de décrire.*

*Preuve.* C'est la traduction, *via* l'isomorphisme entre $U \otimes_{U^e} (U \otimes \Lambda_Q^* V \otimes U)$ et $U \otimes \Lambda_Q^* V$, de la proposition 1.5.2. La différentielle $d$ se calcule à partir de (1) ; dans [27] elle est notée $\widetilde{d}$ et donnée en exemple en 10.2(2). □

## 1.6 Cas particulier de l'espace affine quantique ($r = 0$).

Le premier cas significatif où s'applique cette construction est l'espace affine quantique $\mathcal{O}_Q(k^m)$, qui s'obtient en prenant pour $f$ l'application identiquement nulle, et qui correspond à l'algèbre $A_{m,0}^Q$ (voir 1.3.1 ci-dessus). On rappelle, avec des notations adaptées au présent travail, les résultats obtenus dans ce cas par M. Wambst, et donnés dans la section 6 de [27]. Ces résultats nous serviront à réduire le complexe décrit au paragraphe précédent, pour obtenir un complexe plus petit, calculant encore l'homologie de Hochschild des algèbres $A_{n,r}^\Lambda$ générales.
L'espace affine quantique $\mathcal{O}_Q(k^m)$ paramétré par la matrice $Q = (q_{i,j})$ n'est autre que l'algèbre symétrique $S_Q V$ introduite en remarque après la définition 1.5.1, pour $V = k^m$. Puisque $S_Q V = U_Q(V, f)$ avec $f = 0$, les hypothèses de la proposition 1.5.3 sont vérifiées, et l'homologie de Hochschild de $S_Q V$ est donnée par le complexe suivant :

$$K(S_Q V) = \bigoplus_{\substack{\alpha \in \mathbb{N}^m \\ \gamma \in \{0,1\}^m}} k.v^\alpha \otimes v^\gamma,$$

$$d(v^\alpha \otimes v^\gamma) = \sum_{i=1}^{m} \Omega_Q(\alpha, \gamma; i) v^{\alpha+[i]} \otimes v^{\gamma-[i]}, \tag{5}$$

avec les coefficients :

$$\Omega_Q(\alpha, \gamma; i) = (-1)^{\sum_{k<i} \gamma_k} \left( \prod_{k<i} q_{k,i}^{\gamma_k} \right) \left( \prod_{k>i} q_{k,i}^{\alpha_k} \right) \times \left( 1 - \left( \prod_{k=1}^{m} q_{i,k}^{\alpha_k+\gamma_k} \right) \right) \text{ si } \gamma_i = 1,$$

et

$$\Omega_Q(\alpha, \gamma; i) = 0 \text{ si } \gamma_i = 0; \tag{6}$$

**Proposition 1.6.1** *Le complexe $(K(S_Q V, d))$ ci-dessus calcule l'homologie de Hochschild de $S_Q V$, et il admet une $\mathbb{N}^m$-graduation définie par $deg(v^\alpha \otimes v^\gamma) = \alpha + \gamma$. Par ailleurs, notons $C(Q) = \{\rho \in \mathbb{N}^m \mid \forall i,\ \rho_i = 0 \text{ ou } v_i v^\rho = v^\rho v_i\}$. Alors pour tout $\rho \in \mathbb{N}^m \setminus C(Q)$, le sous-complexe homogène de degré $\rho$ de $(K(S_Q V), d)$ est acyclique.*

*Preuve.* L'algèbre $S_Q V$ vérifie clairement les hypothèses de la proposition 1.5.3, et les formules (2), (3) et (4) donnent bien dans ce cas particulier les formules (5) et (6) pour la différentielle. La graduation découle directement de la formule (5), et l'acyclicité pour $\rho \notin C(Q)$ se montre à l'aide d'une homotopie, décrite dans la démonstration du théorème 6.1 de [27]. □

## 1.7 Cas de l'algèbre $A_{n,r}^{\Lambda}$ pour $r$ quelconque.

Dans ce paragraphe on présente $A_{n,r}^{\Lambda}$ sous forme d'une algèbre enveloppante quantique, et on démontre au lemme 1.7.2 que la condition de la proposition 1.5.2 est vérifiée.

Précisons la matrice associée à la symétrie de Hecke $c$ conformément à la définition 1.5.1 et à la remarque qui la suit. Etant donnée $\Lambda = (\lambda_{i,j})_{1\leq i,j\leq n}$ la matrice introduite en 1.1, on en extrait $\Lambda_r = (\lambda_{i,j})_{1\leq i,j\leq r}$. C'est une matrice multiplicativement antisymétrique de taille $r$. Par ailleurs, pour tous $0 \leq k,t \leq n$ posons

$$\Lambda_{k,t}^{-1} = (\lambda_{i,j}^{-1})_{\substack{1\leq i\leq k \\ 1\leq j\leq t}} \in M_{k,t}(k^*).$$

Définissons par bloc la matrice $Q(\Lambda)$ multiplicativement antisymétrique de taille $n+r$ :

$$Q(\Lambda) = \begin{pmatrix} \Lambda_r & \Lambda_{r,n}^{-1} \\ \Lambda_{n,r}^{-1} & \Lambda \end{pmatrix}. \tag{7}$$

**Proposition 1.7.1** *L'algèbre $A_{n,r}^{\Lambda}$ définie en 1.1 est isomorphe à l'algèbre enveloppante quantique $U_Q(V,f)$, où $Q = Q(\Lambda)$ est définie ci-dessus,*

$$V = k.x_1 \oplus \ldots k.x_r \oplus k.y_1 \oplus \ldots k.y_n,$$

*et $f$ est définie par*

$$f(x_i \otimes x_i) = f(y_j \otimes y_j) = 0 \text{ pour tous } i,j,$$

$$f(x_i \otimes x_j) = f(y_i \otimes y_j) = f(x_i \otimes y_j) = f(y_i \otimes x_j) = 0 \text{ pour } i \neq j,$$

*et*

$$f(x_i \otimes y_i) = -f(y_i \otimes x_i) = 1.$$

*Preuve.* On vérifie aisément que les conditions de la définition 1.5.1 sont satisfaites, et que l'algèbre $U_Q(V,f)$ est isomorphe à $A_{n,r}^{\Lambda}$. □

**Lemme 1.7.2** *Soit $f$ la forme linéaire définie à la proposition 1.7.1. Pour tout $2 \leq k \leq n$ posons $\Pi_k = c_1 \circ \ldots c_{k-1}$, et $\Pi_1 = Id$. De même, pour tout entier $1 \leq k \leq n-1$ posons $\check{\Pi}_k = c_{n-1} \circ \ldots c_k$, et $\check{\Pi}_n = Id$. Alors pour tout $p \geq 2$ l'application $f'$ définie de $V^{\otimes p}$ dans $V^{\otimes p-2}$ par :*

$$f' = \sum_{1\leq i<j\leq p} (-1)^{i+j+1} \left((f \otimes I_{n-2})(I_1 \otimes \Pi_{j-1})\Pi_i - (I_{n-2} \otimes f)(\check{\Pi}_i \otimes I_1)\check{\Pi}_j\right), \tag{8}$$

*est nulle.*

*Preuve.* Il suffit de le démontrer pour des monômes $a_1 \otimes \ldots \otimes a_p$, avec $a_s$ dans l'ensemble $\{x_1, \ldots, x_r, y_1, \ldots, y_n\}$ pour tout $s$. Par définition de $f$ et de $c$, pour $i < j$, le terme d'indice $(i,j)$ dans la somme (8) est nul, sauf peut-être si on a $x_k$ en position $i$ et $y_k$ en $j$, ou l'inverse. Calculons donc ce terme sur

$$X = a_1 \otimes \ldots a_{i-1} \otimes x_k \otimes a_{i+1} \ldots a_{j-1} \otimes y_k \otimes a_{j+1} \ldots a_p,$$

avec $a_s \in \{x_1, \ldots, x_r, y_1, \ldots, y_n\}$. Alors la définition de $c$ implique que

$$(I_1 \otimes \Pi_{j-1})\Pi_i(X) = \pi x_k \otimes y_k \otimes a_1 \otimes \ldots a_{i-1} \otimes a_{i+1} \ldots a_{j-1} \otimes a_{j+1} \ldots a_p,$$

et

$$(\check{\Pi}_i \otimes I_1)\check{\Pi}_j(X) = \check{\pi} a_1 \otimes \ldots a_{i-1} \otimes a_{i+1} \ldots a_{j-1} \otimes a_{j+1} \ldots a_n \otimes x_k \otimes y_k,$$

$$\text{avec} \quad \pi = \prod_{s<i} \lambda_{t(s),k}^{\varepsilon(s)} \times \prod_{s<j} \lambda_{t(s),k}^{-\varepsilon(s)},$$

$$\text{et} \quad \check{\pi} = \prod_{s>j} \lambda_{k,t(s)}^{-\varepsilon(s)} \times \prod_{s>i} \lambda_{k,t(s)}^{\varepsilon(s)} ;$$

où l'on a posé :

$$\begin{array}{lll} t(s) = l \text{ et } \varepsilon(s) = 1 & \text{si } a_s = x_l & \text{pour } s \notin \{i, j\}, \\ t(s) = l \text{ et } \varepsilon(s) = -1 & \text{si } a_s = y_l & \text{pour } s \notin \{i, j\}, \\ t(i) = t(j) = k. & & \end{array}$$

Ainsi :

$$\left((f \otimes I_{n-2})(I_1 \otimes \Pi_{j-1})\Pi_i - (I_{n-2} \otimes f)(\check{\Pi}_i \otimes I_1)\check{\Pi}_j\right)(X) =$$
$$\left(\prod_{i<s<j} \lambda_{t(s),k}^{-\varepsilon(s)} - \prod_{i<s<j} \lambda_{k,t(s)}^{\varepsilon(s)}\right) a_1 \otimes \ldots a_{i-1} \otimes a_{i+1} \ldots a_{j-1} \otimes a_{j+1} \ldots a_n.$$

La matrice $\Lambda$ étant multiplicativement antisymétrique, la parenthèse dans le second terme de cette égalité est nulle. Le même raisonnement s'applique lorsque $y$ est en place $i$ et $x$ en place $j$, d'où le résultat. □

On va donc pouvoir utiliser le complexe de la proposition 1.5.3 pour calculer l'homologie de Hochschild de $A_{n,r}^{\Lambda}$. On conviendra de noter, pour tous $\alpha \in \mathbb{N}^r$, $\beta \in \mathbb{N}^n$, $\gamma \in \{0,1\}^r$, $\delta \in \{0,1\}^n$,

$$x^\alpha y^\beta \otimes x^\gamma y^\delta = x^\alpha y^\beta \otimes x_1^{\gamma_1} \wedge \ldots x_r^{\gamma_r} \wedge y_1^{\delta_1} \wedge \ldots y_n^{\delta_n}.$$

Rappelons que $[i]$ désigne le $i^{\text{ème}}$ élément de la base canonique de $\mathbb{N}^r$ ou $\mathbb{N}^n$. Le degré différentiel d'un monôme $x^\alpha y^\beta \otimes x^\gamma y^\delta$ est égal à $|\gamma + \delta|$. Soit $K = (K(A_{n,r}^{\Lambda}), d)$ le complexe de la proposition 1.5.3. Ce complexe s'écrit :

$$K = \bigoplus_{\substack{\alpha \in \mathbb{N}^r, \beta \in \mathbb{N}^n \\ \gamma \in \{0,1\}^r, \delta \in \{0,1\}^n}} k.x^\alpha y^\beta \otimes x^\gamma y^\delta, \tag{9}$$

où la différentielle $d$ est donnée par la formule :

$$\begin{aligned} & d(x^\alpha y^\beta \otimes x^\gamma y^\delta) = \\ & \sum_{i=1}^{r} \Omega_1(\alpha,\beta,\gamma,\delta;i) x^{\alpha+[i]} y^\beta \otimes x^{\gamma-[i]} y^\delta + \Omega_1'(\alpha,\beta,\gamma,\delta;i) x^\alpha y^{\beta-[i]} \otimes x^{\gamma-[i]} y^\delta \\ & + \sum_{j=1}^{r} \Omega_2(\alpha,\beta,\gamma,\delta;j) x^\alpha y^{\beta+[j]} \otimes x^\gamma y^{\delta-[j]} + \Omega_2'(\alpha,\beta,\gamma,\delta;j) x^{\alpha-[j]} y^\beta \otimes x^\gamma y^{\delta-[j]} \\ & + \sum_{j=r+1}^{n} \Omega_2(\alpha,\beta,\gamma,\delta;j) x^\alpha y^{\beta+[j]} \otimes x^\gamma y^{\delta-[j]}, \end{aligned} \tag{10}$$

avec les coefficients suivants : pour $1 \leq i \leq r$,

$$\begin{aligned} \Omega_1(\alpha,\beta,\gamma,\delta;i) = & \ \epsilon_1(i) \left( \prod_{k<i} \lambda_{k,i}^{\gamma_k} \right) \left( \prod_{k=1}^{n} \lambda_{k,i}^{-\beta_k} \right) \left( \prod_{k=i+1}^{r} \lambda_{k,i}^{\alpha_k} \right) \\ & \times \left( 1 - \left( \prod_{k=1}^{r} \lambda_{i,k}^{\alpha_k+\gamma_k} \right) \left( \prod_{k=1}^{n} \lambda_{i,k}^{-(\beta_k+\delta_k)} \right) \right) \quad \text{si } \gamma_i = 1, \\ \text{et} & \\ \Omega_1(\alpha,\beta,\gamma,\delta;i) = & \ 0 \quad \text{sinon;} \end{aligned} \tag{11}$$

pour $1 \leq i \leq r$,

$$\begin{aligned} \Omega_1'(\alpha,\beta,\gamma,\delta;i) = & \ -\epsilon_1(i)\beta_i \left( \prod_{k<i} \lambda_{k,i}^{\gamma_k} \right) \left( \prod_{k=i+1}^{n} \lambda_{k,i}^{-\beta_k} \right) \quad \text{si } \gamma_i = 1, \\ \text{et} & \\ \Omega_1'(\alpha,\beta,\gamma,\delta;i) = & \ 0 \quad \text{sinon,} \end{aligned} \tag{12}$$

avec $\epsilon_1(i) = (-1)^{\sum_{k<i} \gamma_k}$.
Pour $1 \leq j \leq n$,

$$\begin{aligned} \Omega_2(\alpha,\beta,\gamma,\delta;j) = & \ \epsilon_2(j) \left( \prod_{k<j} \lambda_{k,j}^{\delta_k} \right) \left( \prod_{k=1}^{r} \lambda_{k,j}^{-\gamma_k} \right) \left( \prod_{k=j+1}^{n} \lambda_{k,j}^{\beta_k} \right) \\ & \times \left( 1 - \left( \prod_{k=1}^{r} \lambda_{j,k}^{-(\alpha_k+\gamma_k)} \right) \left( \prod_{k=1}^{n} \lambda_{j,k}^{\beta_k+\delta_k} \right) \right) \quad \text{si } \delta_j = 1, \\ \text{et} & \\ \Omega_2(\alpha,\beta,\gamma,\delta;j) = & \ 0 \quad \text{sinon;} \end{aligned} \tag{13}$$

et pour $1 \leq j \leq r$,

$$\Omega'_2(\alpha,\beta,\gamma,\delta;j) = \epsilon_2(j)\alpha_j \left(\prod_{k=j+1}^{n} \lambda_{j,k}^{\delta_k}\right)\left(\prod_{k<j} \lambda_{j,k}^{-\alpha_k}\right) \quad \text{si } \delta_j = 1, \tag{14}$$

et

$$\Omega'_2(\alpha,\beta,\gamma,\delta;j) = 0 \quad \text{sinon},$$

avec $\epsilon_2(j) = (-1)^{|\gamma| + \sum_{k<j} \delta_k}$.

**Proposition 1.7.3** *L'homologie de Hochschild de $A^{\Lambda}_{n,r}$ est l'homologie du complexe décrit ci-dessus :*

$$HH_*(A^{\Lambda}_{n,r}) = H_*(K,d).$$

*Preuve.* Par la proposition 1.7.1 et le lemme 1.7.2, on peut appliquer la proposition 1.5.3 à l'algèbre $A^{\Lambda}_{n,r}$. □

**Corollaire 1.7.4** *L'homologie de Hochschild de $A^{\Lambda}_{n,r}$ est nulle en degré strictement supérieur à $n + r$.*

*Preuve.* Par définition les espaces $K_*$ sont nuls pour $* > n + r$. □

# 2 Réduction du complexe $K$.

Nous montrons dans cette section que le complexe de Koszul quantique $(K,d)$ de l'algèbre $A^{\Lambda}_{n,r}$ contient un sous-complexe $(K_C,d)$ tel que le complexe quotient $K/K_C$ est acyclique. Il en résulte que les complexes $K$ et $K_C$ ont la même homologie. Considérons donc le complexe $(K,d)$ défini par (9) et (10).

## 2.1 Graduation de $K$.

Définissons le degré total d'un monôme de $K$ par :

$$\deg(x^{\alpha}y^{\beta} \otimes x^{\gamma}y^{\delta}) = (\alpha+\gamma, \beta+\delta) \in \mathbb{N}^{r+n}.$$

Pour tout $\rho \in \mathbb{N}^{r+n}$, notons $K_\rho$ la partie homogène de degré total $\rho$ de $K$. On constate que la différentielle $d$ ne respecte pas cette graduation de $K$, mais qu'en général on a le résultat suivant :

**Lemme 2.1.1** $d(K_\rho) \subset K_\rho \oplus \bigoplus_{i=1}^{r} K_{\rho-[i]-[r+i]}$.

*Preuve.* Ceci découle directement de la formule (10). □

Afin d'utiliser les résultats d'acyclicité de la proposition 1.6.1, rappelons la matrice $Q(\Lambda)$ définie en (7) :

$$Q(\Lambda) = \begin{pmatrix} \Lambda_r & \Lambda_{r,n}^{-1} \\ \Lambda_{n,r}^{-1} & \Lambda \end{pmatrix}.$$

Notons $(\widetilde{\lambda}_{k,i})_{1\le k,i\le n+r}$ les élements de $Q(\Lambda)$, de sorte que $\widetilde{\lambda}_{k,i} = Q(\Lambda)_{k,i}$. Conformément à la notation introduite dans la proposition 1.6.1, posons alors la définition suivante :

**Définition 2.1.2** *Soit $C$ la partie de $\mathbb{N}^{r+n}$ définie par :*

$$C = C(Q(\Lambda)) = \{\rho \in \mathbb{N}^{r+n} \mid \forall i,\ \rho_i = 0 \ ou \ \prod_{k=1}^{r+n} \widetilde{\lambda}_{k,i}^{\rho_k} = 1\}.$$

**Lemme 2.1.3** *1. L'ensemble $C$ est constitué des $\rho \in \mathbb{N}^{r+n}$ tels que :*

– *pour tout $i \leq r$, on a $\prod_{k=1}^{r+n} \widetilde{\lambda}_{k,i}^{\rho_k} = 1$ ou $\rho_i = \rho_{i+r} = 0$,*

– *pour tout $i \geq 2r+1$, $\prod_{k=1}^{r+n} \widetilde{\lambda}_{k,i}^{\rho_k} = 1$ ou $\rho_i = 0$.*

*2. Si $\rho = (\alpha+\gamma, \beta+\delta) \in C$, alors $\Omega_1(\alpha,\beta,\gamma,\delta;i) = \Omega_2(\alpha,\beta,\gamma,\delta;j) = 0$ pour tous $i, j$.*

*Preuve.* 1. Ceci découle du fait que dans $Q(\Lambda)$, pour tout $i \leq r$, et pour tout $1 \leq k \leq r+n$, on a $\widetilde{\lambda}_{i,k} = (\widetilde{\lambda}_{i+r,k})^{-1}$.
2. C'est une conséquence directe de la définition de $C$ et des formules (11) et (13).
□
On définit alors les deux sous-espaces vectoriels suivants de $K$ :

$$K_C = \bigoplus_{\rho \in C} K_\rho, \quad K'_C = \bigoplus_{\rho \notin C} K_\rho.$$

**Lemme 2.1.4** *La différentielle $d$ vérifie $d(K_C) \subset K_C$.*

*Preuve.* Soit $\rho \in C$, et $x^\alpha y^\beta \otimes x^\gamma y^\delta$ un monôme de degré $\rho$. Alors par le point 2 du lemme 2.1.3, dans $d(x^\alpha y^\beta \otimes x^\gamma y^\delta)$ les termes $\Omega_1(\alpha,\beta,\gamma,\delta;i)$ et $\Omega_2(\alpha,\beta,\gamma,\delta;j)$ sont nuls pour tous $i, j$, et les autres termes sont de degré $\rho - [i] - [r+i]$ si $\gamma_i = 1$ ou $\delta_i = 1$, dont on vérifie facilement qu'ils sont encore dans $C$. □

On note encore $d$ la restriction de la différentielle de $K$ à $K_C$. Notons $(\overline{K}, \overline{d})$ le complexe quotient. De la suite exacte courte de complexe

$$0 \longrightarrow K_C \xrightarrow{\iota} K \xrightarrow{p} \overline{K} \longrightarrow 0,$$

où $\iota$ et $p$ sont les applications canoniques, on déduit la suite exacte longue d'homologie (voir par exemple [19], th. 4.1) :

$$\cdots \longrightarrow H_{*+1}(\overline{K}) \longrightarrow H_*(K_C) \longrightarrow H_*(K) \longrightarrow H_*(\overline{K}) \longrightarrow \cdots \qquad (15)$$

Nous allons montrer que le complexe quotient $(\overline{K}, \overline{d})$ est acyclique, ce qui grâce à la suite exacte longue (15) montrera que les complexes $(K, d)$ et $(K_C, d)$ ont la même homologie.

## 2.2 Acyclicité de $(\overline{K}, \overline{d})$.

Montrons que le complexe gradué de $\overline{K}$ pour une certaine filtration est acyclique. Le $k$-espace vectoriel $\overline{K}$ est isomorphe à $K'_C$. Définissons alors une $\mathbb{N}$-filtration $\mathcal{F}$ de $(\overline{K}, \overline{d})$ en posant pour tout $t \in \mathbb{N}$ :

$$\mathcal{F}^t(\overline{K}) = \bigoplus_{\substack{|\alpha| + |\beta| + |\gamma| + |\delta| \leq t \\ (\alpha+\gamma, \beta+\delta) \notin C}} k.\overline{x}^\alpha \overline{y}^\beta \otimes \overline{x}^\gamma \overline{y}^\delta.$$

Il découle des définitions et de la formule (10) que ceci définit une filtration de complexe pour $\overline{K}$. Notons $\widetilde{K} = \mathrm{gr}(\overline{K})$ le complexe gradué, et $\widetilde{d} = \mathrm{gr}(\overline{d})$ sa différentielle. Notons également $\widetilde{x}^\alpha \widetilde{y}^\beta \otimes \widetilde{x}^\gamma \widetilde{y}^\delta = \mathrm{gr}(\overline{x}^\alpha \overline{y}^\beta \otimes \overline{x}^\gamma \overline{y}^\delta)$. On a clairement :

$$\widetilde{K} = \bigoplus_{(\alpha+\gamma, \beta+\delta) \notin C} k.\widetilde{x}^\alpha \widetilde{y}^\beta \otimes \widetilde{x}^\gamma \widetilde{y}^\delta,$$

$$\begin{aligned} \widetilde{d}(\widetilde{x}^\alpha \widetilde{y}^\beta \otimes \widetilde{x}^\gamma \widetilde{y}^\delta) = & \sum_{i=1}^{r} \Omega_1(\alpha, \beta, \gamma, \delta; i) \widetilde{x}^{\alpha+[i]} \widetilde{y}^\beta \otimes \widetilde{x}^{\gamma-[i]} \widetilde{y}^\delta \\ & + \sum_{j=1}^{n} \Omega_2(\alpha, \beta, \gamma, \delta; j) \widetilde{x}^\alpha \widetilde{y}^{\beta+[j]} \otimes \widetilde{x}^\gamma \widetilde{y}^{\delta-[j]}. \end{aligned} \tag{16}$$

On introduit alors une $\mathbb{N}^{r+n}$ graduation de l'espace $\widetilde{K}$ de la façon suivante : pour un multi-entier $\rho$ de $\mathbb{N}^{r+n}$ n'appartenant pas à $C$, on pose

$$\widetilde{K_\rho} = \bigoplus_{(\alpha+\gamma, \beta+\delta) = \rho} k.\widetilde{x}^\alpha \widetilde{y}^\beta \otimes \widetilde{x}^\gamma \widetilde{y}^\delta,$$

de sorte que

$$\widetilde{K} = \bigoplus_{\rho \notin C} \widetilde{K_\rho}.$$

La formule (16) montre que la différentielle $\widetilde{d}$ respecte cette graduation, et $(\widetilde{K_\rho}, \widetilde{d})$ est un sous-complexe de $\widetilde{K}$.

**Lemme 2.2.1** *Chacune des composantes homogènes $(\widetilde{K_\rho}, \widetilde{d})$ est acyclique.*

*Preuve.* En effet, soient $V = k^{r+n}$, et $(K(S_{Q(\Lambda)}V), d_{Q(\Lambda)})$ le complexe associé à l'algèbre symétrique $S_{Q(\Lambda)}V$ au paragraphe 1.6, en prenant pour matrice $Q$ la matrice $Q(\Lambda)$. Alors pour tout $\rho \notin C$, le complexe $(\widetilde{K_\rho}, \widetilde{d_\rho})$ est isomorphe à la partie homogène de degré $\rho$ du complexe $(K(S_{Q(\Lambda)}V), d_{Q(\Lambda)})$. Par ailleurs, par la définition 2.1.2 on a $C = C(Q(\Lambda))$, et on conclut grâce à la proposition 1.6.1. □

**Corollaire 2.2.2** *Les complexes $(K, d)$ et $(K_C, d)$ ont même homologie.*

*Preuve.* Il découle du lemme précédent que le complexe $(\widetilde{K}, \widetilde{d}) = \mathrm{gr}(\overline{K}, \overline{d})$ est acyclique. La filtration $\mathcal{F}$ étant croissante, exhaustive et bornée inférieurement, on en déduit l'acyclicité du complexe $(\overline{K}, \overline{d})$ (ce résultat, classique en algèbre homologique, peut se démontrer en raisonnant par récurrence sur le degré filtrant d'un élément de $\overline{K}$). La suite exacte longue (15) donne alors pour tout $* \in \mathbb{N}$ une suite exacte courte

$$0 \longrightarrow H_*(K_C) \longrightarrow H_*(K) \longrightarrow 0,$$

ce qui termine la preuve. □

### 2.3 Homologie de $A_{n,r}^{\Lambda}$.

Par construction, le petit complexe $(K_C, d)$ se présente de la façon suivante :

$$K_C = \bigoplus_{(\alpha+\gamma,\beta+\delta)\in C} k.x^{\alpha}y^{\beta} \otimes x^{\gamma}y^{\delta}, \tag{17}$$

$$\begin{aligned} d(x^{\alpha}y^{\beta} \otimes x^{\gamma}y^{\delta}) = \sum_{i=1}^{r} \Omega'_1(\alpha,\beta,\gamma,\delta;i)x^{\alpha}y^{\beta-[i]} \otimes x^{\gamma-[i]}y^{\delta} \\ +\sum_{j=1}^{r} \Omega'_2(\alpha,\beta,\gamma,\delta;j)x^{\alpha-[j]}y^{\beta} \otimes x^{\gamma}y^{\delta-[j]}. \end{aligned} \tag{18}$$

Enonçons le résultat principal de cette section, qui permettra dans les trois prochaines sections de faire des calculs explicites.

**Théorème 2.3.1** *Soient $r \leq n$ deux entiers, et $\Lambda \in M_n(k^*)$ une matrice multiplicativement antisymétrique. Alors l'homologie de Hochschild de l'algèbre $A_{n,r}^{\Lambda}$ définie en 1.1 est l'homologie du petit complexe $(K_C, d)$ décrit par les équations (17) et (18), où l'ensemble $C$ est défini en 2.1.2.*

*Preuve.* Ceci découle directement de la proposition 1.7.3 et du corollaire 2.2.2. □

Nous allons maintenant appliquer le complexe $(K_C, d)$ au calcul de l'homologie de Hochschild des algèbres $A_{n,r}^{\Lambda}$ dans certains cas particuliers : le cas "semi-classique", le cas "libre" et le cas "mixte minimal" ($r = 1, n = 2$).

## 3 Application au cas semi-classique ($n = r$).

Nous considérons dans cette section l'algèbre $A_{n,n}^{\Lambda}$ définie en 1.1, avec $r = n$. Comme on l'a déjà remarqué en 1.3.3 il s'agit alors d'un cas particulier d'algèbre de Weyl quantique $A_n^{(1),\Lambda}(k)$, algèbre d'opérateurs différentiels sur un espace affine quantique paramétré par $\Lambda$.

Nous allons montrer que dans cette situation l'homologie du complexe $K_C$ est celle de l'algèbre de Weyl classique.

## 3.1 Rappels sur l'algèbre de Weyl classique $A_n(k)$.

L'algèbre de Weyl $A_n(k)$ est engendrée par $2n$ générateurs $p_1, \ldots, p_n, q_1, \ldots, q_n$ soumis aux relations :

$$[p_i, p_j] = [q_i, q_j] = [p_i, q_j] = 0 \text{ pour } i \neq j,\ [p_i, q_i] = 1.$$

C'est un cas particulier de l'algèbre $A_{n,n}^{\Lambda}$, où tous les $\lambda_{i,j}$ sont égaux à 1. La proposition 1.7.3 s'applique donc, et le complexe $(K, d)$ décrit dans la section 1.7, calcule l'homologie de Hochschild de $A_n(k)$. Dans ce cas ce complexe, noté $(K_W, d_W)$, devient :

$$K_W = \bigoplus_{\substack{\alpha,\, \beta \,\in\, \mathbf{N}^n \\ \gamma,\, \delta \,\in\, \{0,1\}^n}} p^\alpha q^\beta \otimes p^\gamma q^\delta.$$

$$\begin{aligned} d_W(p^\alpha q^\beta \otimes p^\gamma q^\delta) = &\sum_{\gamma_i=1} \epsilon_1(i)(-\beta_i) p^\alpha q^{\beta-[i]} \otimes p^{\gamma-[i]} q^\delta \\ &+ \sum_{\delta_j=1} \epsilon_2(j)(-\alpha_j) p^{\alpha-[j]} q^\beta \otimes p^\gamma q^{\delta-[j]}. \end{aligned} \tag{19}$$

Ceci montre que le complexe $(K_W, d_W)$ est exactement le complexe de Koszul usuel de l'algèbre de Weyl $A_n(k)$.
Rappelons l'homologie de Hochschild de l'algèbre de Weyl (voir par exemple [26], [17], [21]) :

$$\begin{aligned} HH_*(A_n(k)) &= 0 \ \text{ si } * \neq 2n, \\ HH_{2n}(A_n(k)) &= k. \end{aligned}$$

## 3.2 Lien entre les complexes $K_C$ et $K_W$.

Définissons une application linéaire $f$ de $K_C$ dans $K_W$ par :

$$f(x^\alpha y^\beta \otimes x^\gamma y^\delta) = R(\alpha, \beta, \gamma, \delta) p^\alpha q^\beta \otimes p^\gamma q^\delta, \tag{20}$$

où

$$R(\alpha, \beta, \gamma, \delta) = \left( \prod_{i,\ \gamma_i=0} \prod_{k<i} \lambda_{k,i}^{-\gamma_k} \prod_{k>i} \lambda_{k,i}^{\beta_k} \right) \times \left( \prod_{j,\ \delta_j=0} \prod_{k>j} \lambda_{j,k}^{-\delta_k} \prod_{k<j} \lambda_{j,k}^{\alpha_k} \right) \in k^*. \tag{21}$$

**Proposition 3.2.1** *L'application $f : (K_C, d) \to (K_W, d_W)$ est un morphisme de complexes.*

*Preuve.* On doit vérifier que pour tous $(\alpha, \beta) \in \mathbb{N}^{r+n}$, et $(\gamma, \delta) \in \{0,1\}^{r+n}$ tels que $(\alpha + \gamma, \beta + \delta) \in C$, on a :

$$f \circ d(x^\alpha y^\beta \otimes x^\gamma y^\delta) = d_W \circ f(x^\alpha y^\beta \otimes x^\gamma y^\delta). \tag{22}$$

Les formules (20), (19) et (18) permettent de calculer :

$$f \circ d(x^\alpha y^\beta \otimes x^\gamma y^\delta) =$$
$$\sum_{i=1}^{n} \Omega'_1(\alpha,\beta,\gamma,\delta,i) R(\alpha,\beta-[i],\gamma-[i],\delta) p^\alpha q^{\beta-[i]} \otimes p^{\gamma-[i]} q^\delta$$
$$+\sum_{j=1}^{n} \Omega'_2(\alpha,\beta,\gamma,\delta,j) R(\alpha-[j],\beta,\gamma,\delta-[j]) p^{\alpha-[j]} q^\beta \otimes p^\gamma q^{\delta-[j]}.$$

Par ailleurs, on a :

$$\begin{aligned} d_W \circ f(x^\alpha y^\beta \otimes x^\gamma y^\delta) &= \sum_{\gamma_i=1} R(\alpha,\beta,\gamma,\delta)\epsilon_1(i)(-\beta_i) p^\alpha q^{\beta-[i]} \otimes p^{\gamma-[i]} q^\delta \\ &\quad + \sum_{\delta_j=1} R(\alpha,\beta,\gamma,\delta)\epsilon_2(j)(-\alpha_j) p^{\alpha-[j]} q^\beta \otimes p^\gamma q^{\delta-[j]}. \end{aligned}$$

L'équation (22) est satisfaite si et seulement si pour tout $i$ tel que $\gamma_i = 1$ on a :

$$\Omega'_1(\alpha,\beta,\gamma,\delta,i) R(\alpha,\beta-[i],\gamma-[i],\delta) = R(\alpha,\beta,\gamma,\delta)\epsilon_1(i)(-\beta_i),$$

et pour tout $j$ tel que $\delta_j = 1$ :

$$\Omega'_2(\alpha,\beta,\gamma,\delta,j) R(\alpha-[j],\beta,\gamma,\delta-[j]) = R(\alpha,\beta,\gamma,\delta)\epsilon_2(j)(-\alpha_j),$$

ce qui se vérifie directement à partir des formules (12), (14) et (21) définissant les scalaires

$$\Omega'_1(\alpha,\beta,\gamma,\delta;i),\ \Omega'_2(\alpha,\beta,\gamma,\delta;j) \text{ et } R(\alpha,\beta,\gamma,\delta).$$

□

On définit ensuite une application linéaire $g$ de $K_W$ dans $K_C$ par :

$g(p^\alpha q^\beta \otimes p^\gamma q^\delta) = R(\alpha,\beta,\gamma,\delta)^{-1} x^\alpha y^\beta \otimes x^\gamma y^\delta$ si $(\alpha+\gamma, \beta+\delta) \in C$,
et
$g(p^\alpha q^\beta \otimes p^\gamma q^\delta) = 0$ sinon.

Des calculs analogues à ceux de la proposition 3.2.1 montrent que $g$ est également un morphisme de complexes.

**Proposition 3.2.2** *Dans le diagramme ci-dessous, $f$ et $g$ sont des morphismes de complexes tels que $g \circ f = Id_{K_C}$.*

$$(K_C, d) \underset{g}{\overset{f}{\rightleftharpoons}} (K_W, d_W).$$

*En particulier, il existe une injection linéaire de $H_*(K_C)$ dans $H_*(K_W)$.*

*Preuve.* Puisque $f$ et $g$ sont des morphismes de complexes différentiels, ils induisent des applications linéaires $H(f)$ et $H(g)$ sur les homologies. En outre par fonctorialité de la relation $g \circ f = \mathrm{Id}_{K_C}$ on déduit $H(g) \circ H(f) = \mathrm{Id}_{H(K_C)}$. Donc $H(f)$ est l'injection annoncée. □

**Corollaire 3.2.3**

$$\begin{array}{ll} H_*(K_C,d) & = 0 \quad si * \neq 2n, \\ H_{2n}(K_C,d) & = k. \end{array}$$

*Preuve.* L'homologie de l'algèbre de Weyl $A_n(k)$ est nulle sauf en degré $2n$, où elle vaut $k$ (voir ci-dessus la fin du paragraphe 3.1). D'après la proposition 3.2.2 il suffit de vérifier que $(K_C,d)$ a une homologie non nulle en degré $2n$. Or il est clair que le multi-entier $\rho_o = (1,\ldots,1)$ appartient à $C$. Considérons alors le monôme $1 \otimes x_1 \wedge \ldots \wedge y_n$ : il est de degré $\rho_o$, son degré différentiel est égal à $2n$, et $d(1 \otimes x_1 \wedge \ldots \wedge y_n) = 0$. La forme du complexe $(K_C,d)$ implique que

$$H_{2n}(K_C) = \operatorname{Ker} d_{2n} / \operatorname{Im} d_{2n+1} = \operatorname{Ker} d_{2n}.$$

Ce noyau contient au moins $1 \otimes x_1 \wedge \ldots \wedge y_n$, il n'est donc pas réduit à 0, et le résultat est démontré. □

Enonçons pour conclure le résultat principal de cette section.

**Théorème 3.2.4 (Homologie des algèbres semi-classiques $A^{\Lambda}_{n,n}$)** *Soient $n$ un entier non nul, et $\Lambda \in M_n(k^*)$ une matrice multiplicativement antisymétrique. Alors l'homologie de Hochschild de l'algèbre $A^{\Lambda}_{n,n}$ définie en 1.1 avec $r = n$ est donnée par :*

$$\begin{array}{ll} HH_*(A^{\Lambda}_{n,n}) & = 0 \quad si * \neq 2n, \\ HH_{2n}(A^{\Lambda}_{n,n}) & = k. \end{array}$$

*Preuve.* Ceci découle du théorème 2.3.1 et du corollaire 3.2.3. □

# 4 Application au cas libre.

Cette section est consacrée à la situation suivante. On se donne une famille $(\lambda_{i,j})_{1\leq i<j\leq n}$ dans $k^*$ formant un groupe abélien libre de rang maximal, c'est-à-dire de rang $n(n-1)/2$. Ce sont par exemple les hypothèses prises par V.Artamonov pour définir les “general quantum functions” dans [3].

**Remarque :** pour $n = 2$, le seul paramètre à considérer est $\lambda_{1,2}$. L'hypothèse de liberté équivaut alors à dire que $\lambda_{1,2}$ n'est pas racine de l'unité dans $k^*$.

Pour un $r \leq n$ fixé, considérons l'algèbre $A^{\Lambda}_{n,r}$ définie en 1.1 à partir de la matrice $\Lambda = (\lambda_{i,j}) \in M_n(k^*)$, où l'on a posé $\lambda_{j,i} = \lambda_{i,j}^{-1}$ pour $j > i$ et $\lambda_{i,i} = 1$. Cette matrice est multiplicativement antisymétrique par construction. On définit alors la matrice $Q(\Lambda)$ conformément à la formule (7). Par ailleurs, on note $\Lambda_r$ la matrice extraite de $\Lambda$ en prenant les $r$ premières lignes et colonnes, de sorte que $\Lambda_r = (\lambda_{i,j})_{1\leq i,j\leq r}$ est multiplicativement antisymétrique de taille $r$. Conformément aux constructions faites dans la section précédente, on pose $Q(\Lambda_r) = \begin{pmatrix} \Lambda_r & \Lambda_r^{-1} \\ \Lambda_r^{-1} & \Lambda_r \end{pmatrix}$ avec les conventions du paragraphe 2.1, et on définit une partie $C_r = C(Q(\Lambda_r))$ de $\mathbb{N}^{2r}$ en suivant la définition 2.1.2. Soit alors $(K_{C_r}, d_r)$ le petit complexe qui d'après le théorème 2.3.1 calcule l'homologie de Hochschild de $A^{\Lambda_r}_{r,r}$.

## 4.1 Graduation du complexe $K_C$.

Traduisons tout d'abord la caractérisation de $C$ dans le cas libre.

**Lemme 4.1.1** *Supposons que les $\lambda_{i,j}$ forment un groupe abélien libre de rang maximal, alors l'ensemble $C = C(Q(\Lambda))$ défini en 2.1.2 pour la matrice $\Lambda$ introduite ci-dessus est l'ensemble des $\rho \in \mathbb{N}^{r+n}$ satisfaisant les deux assertions suivantes :*

1. *pour tout $i \leq r$, un des deux points suivants est vérifié :*
   - *soit pour tout $k \neq i$, on a $\rho_k - \rho_{k+r} = 0$ si $k \leq r$ et $\rho_k = 0$ si $k > 2r$,*
   - *soit $\rho_i = \rho_{i+r} = 0$ ;*
2. *pour tout $i \geq 2r+1$, un des deux points suivants est vérifié :*
   - *soit pour tout $k \neq i$, on a $\rho_k - \rho_{k+r} = 0$ si $k \leq r$ et $\rho_k = 0$ si $k > 2r$,*
   - *soit $\rho_i = 0$.*

*Preuve.* Les $\lambda_{i,j}$ étant indépendants dans $k^*$, chaque égalité apparaissant dans la caractérisation de $C$ donnée au lemme 2.1.3 se traduit par la nullité des exposants en chaque $\lambda_{i,j}$. □

Le résultat suivant, indépendant de l'hypothèse de liberté, nous permettra dans le cas libre de calculer l'homologie de $A^{\Lambda}_{n,r}$. C'est pourquoi nous le présentons ici.

**Lemme 4.1.2** *On appelle degré purement quantique d'un monôme $x^\alpha y^\beta \otimes x^\gamma y^\delta$ le multi-entier $(\beta_{r+1} + \delta_{r+1}, \dots, \beta_n + \delta_n) \in \mathbb{N}^{n-r}$. Il définit une $\mathbb{N}^{n-r}$ graduation du complexe $K_C$. Pour tout $\nu \in \mathbb{N}^{n-r}$ on note $K_{C,\nu}$ la partie homogène de $K_C$ de degré purement quantique $\nu$. Alors $K_{C,\nu}$ est un sous-complexe de $K_C$, et l'on a :*

$$K_C = \bigoplus_{\nu \in \mathbb{N}^{n-r}} K_{C,\nu}, \ \text{et} \ H_*(K_C) = \bigoplus_{\nu \in \mathbb{N}^{n-r}} H_*(K_{C,\nu}).$$

*Preuve.* Le degré purement quantique étant défini comme dans le lemme, la formule (18) montre que la différentielle $d$ respecte cette graduation. Le reste en découle immédiatement. □

Montrons que la forme particulière de l'ensemble $C$ donnée par le lemme 4.1.1 permet de calculer facilement l'homologie de $A^{\Lambda}_{n,r}$ dans le cas libre.

## 4.2 Homologie de $K_{C,\nu}$ pour $\nu = (0, \dots, 0)$.

Pour $\nu = (0, \dots, 0)$, on note $K_{C,(0)}$ le complexe $K_{C,\nu}$.

**Proposition 4.2.1** *Le complexe $K_{C,(0)}$ est isomorphe (en tant que complexe différentiel) au complexe $K_{C_r}$ associé à l'algèbre $A^{\Lambda_r}_{r,r}$.*

*Preuve.* D'après le lemme 4.1.1 il est clair que pour un élément $\rho \in \mathbb{N}^n$ tel que $\rho_{2r+1} = \ldots = \rho_{r+n} = 0$, on a équivalence entre $\rho \in C$ et $(\rho_1, \dots, \rho_{2r}) \in C_r$. Ainsi, par définition de la $\mathbb{N}^{n-r}$-graduation de $K_C$ donnée au lemme 4.1.2, les $k$-espaces vectoriels $K_{C,(0)}$ et $K_{C_r}$ sont naturellement isomorphes, et la formule (18) montre que l'isomorphisme respecte les différentielles correspondantes. □

**Corollaire 4.2.2** *L'homologie de $K_{C,(0)}$ est donnée par :*

$$\begin{aligned} HH_*(K_{C,(0)}) &= 0 \text{ si } * \neq 2r, \\ HH_{2r}(K_{C,(0)}) &= k. \end{aligned}$$

*Preuve.* Ceci découle directement de la proposition précédente et du corollaire 3.2.3. □

## 4.3 Homologie de $K_{C,\nu}$ pour $\nu \neq (0, \dots, 0)$.

Le lemme 4.1.1 permet d'obtenir la propriété suivante pour l'ensemble $C$.

**Lemme 4.3.1** *Avec les hypothèses du lemme 4.1.1, soit $\rho = (\rho_1, \dots, \rho_n) \in C$. Supposons qu'il existe un indice $i > 2r$ tel que $\rho_i \neq 0$. Alors $\rho = \rho_i[i]$.*

*Preuve.* Soit $i > 2r$ tel que $\rho_i \neq 0$. Alors, $\rho \in C$ vérifie les assertions 1 et 2 du lemme 4.1.1, et donc pour tout indice $k \leq r$ on a $\rho_k = \rho_{r+k} 0$, et pour tout indice $k > 2r$ distinct de $i$ on a $\rho_k = 0$. □

**Corollaire 4.3.2** *Avec les hypothèses précédentes, pour le sous-complexe de degré purement quantique non nul $\nu \in \mathbb{N}^{n-r} \setminus \{(0, \dots, 0)\}$, on a l'alternative suivante :*

- *soit $K_{C,\nu} = 0$,*
- *soit il existe $p \in \mathbb{N}^*$ et $i > 2r$ tels que $\nu = p.[i - 2r]$, et dans ce cas $d_\nu = 0$, et l'on a :*

$$H_*(K_{C,\nu}, d_\nu) = K_{C,\nu} = k.y_i^p \otimes 1 \;\oplus\; k.y_i^{p-1} \otimes y_i.$$

*Preuve.* Supposons $K_{C,\nu} \neq 0$, et soit $x^\alpha y^\beta \otimes x^\gamma y^\delta$ un monôme appartenant à $K_{C,\nu}$. Alors le degré total $\rho = (\alpha + \gamma, \beta + \gamma)$ de ce monôme vérifie les hypothèses du lemme 4.3.1, donc il existe $i > 2r$ et $p = \rho_i \neq 0$ tels que $\rho = p[i]$, d'où découlent la forme de $\nu$ et de $K_{C,\nu}$. La formule (18) montre alors que la différentielle $d$ est nulle sur la composante homogène $K_{C,\nu}$. □

## 4.4 Homologie des algèbres $A_{n,r}^\Lambda$ dans le cas libre.

**Théorème 4.4.1** *Soit $(\lambda_{i,j})_{1 \leq i < j \leq n}$ une famille d'éléments de $k^*$ formant un groupe abélien libre de rang $n(n-1)/2$. Soit $r \leq n$ un entier, et $\Lambda = (\lambda_{i,j})$ la matrice multiplicativement antisymétrique formée en posant $\lambda_{j,i} = \lambda_{i,j}^{-1}$ et $\lambda_{i,i} = 1$. Alors l'homologie de Hochschild de l'algèbre $A_{n,r}^\Lambda$ définie en 1.1 est la suivante :*
• $1^{er}$ *cas : $r = 0$. Alors :*

$$\begin{aligned} HH_0(A_{n,0}^\Lambda) &= k \;\oplus\; \bigoplus_{i=1}^{n} y_i k[y_i] \otimes 1, \\ HH_1(A_{n,0}^\Lambda) &= \bigoplus_{i=1}^{n} k[y_i] \otimes y_i, \\ HH_*(A_{n,0}^\Lambda) &= 0 \text{ si } * \geq 2; \end{aligned}$$

• $2^{ème}$ *cas : $r \geq 1$. Alors :*

$$\begin{aligned}
HH_0(A^\Lambda_{n,r}) &= \bigoplus_{i=r+1}^{n} y_i k[y_i] \otimes 1, \\
HH_1(A^\Lambda_{n,r}) &= \bigoplus_{i=r+1}^{n} k[y_i] \otimes y_i, \\
HH_{2r}(A^\Lambda_{n,r}) &= k, \\
HH_*(A^\Lambda_{n,r}) &= 0 \quad si \quad * \notin \{0, 1, 2r\}.
\end{aligned}$$

*Preuve.* Par le théorème 2.3.1 et le lemme 4.1.2 l'homologie de Hochschild de $A^\Lambda_{n,r}$ est la somme directe des homologies des complexes $K_{C,\nu}$. Le corollaire 4.3.2 montre que si $\nu \neq (0, \dots, 0)$ alors le complexe $K_{C,\nu}$, est non nul si et seulement si $\nu = p.[i]$, avec $p \in \mathbb{N}^*$, et $i > 2r$, et que dans ce cas il a pour homologie :

$$H_0(K_{C,\nu}) = k.y_i^p \otimes 1,\ H_1(K_{C,\nu}) = k.y_i^{p-1} \otimes y_i,\ H_*(K_{C,\nu}) = 0 \ \text{ si } \ * \geq 2.$$

Ceci, ajouté au résultat du corollaire 4.2.2, termine la preuve. □

**Remarques :** • Dans le cas où $r = n$, les homologies en degré 0 et 1 sont nulles, et on retrouve les résultats du cas semi-classique donnés dans le théorème 3.2.4.
• Dans le cas où $r = 0$, l'algèbre $A^\Lambda_{n,r}$ est un espace affine quantique, et le résultat est un cas particulier du corollaire 6.2 de [27].

# 5 Application au cas mixte minimal ($n = 2,\ r = 1$).

C'est le plus petit cas "mixte" que l'on puisse construire. On note $z$ pour $y_2$, et $\lambda$ pour $\lambda_{2,1}$, de sorte que $A = A^\Lambda_{2,1}$ est engendré par $x, y, z$ avec les relations

$$xy - yx = 1,\ xz = \lambda zx,\ yz = \lambda^{-1} zy.$$

On appelle cette algèbre l'algèbre mixte minimale. On peut facilement vérifier que le centre de cette algèbre est $k$ si $\lambda$ n'est pas racine de l'unité, et $k[z^n]$ si $\lambda$ est d'ordre $n$ dans $k^*$ (voir [25]).
Par construction, le complexe $(K_C, d)$ s'écrit :

$$K_C = \bigoplus_{\alpha+\gamma \in C} x^{\alpha_1} y^{\alpha_2} z^{\alpha_3} \otimes x^{\gamma_1} \wedge y^{\gamma_2} \wedge z^{\gamma_3},$$

$$\begin{aligned}
d(x^{\alpha_1} y^{\alpha_2} z^{\alpha_3} \otimes x^{\gamma_1} \wedge y^{\gamma_2} \wedge z^{\gamma_3}) = \Omega'_1 x^{\alpha_1} y^{\alpha_2 - 1} z^{\alpha_3} \otimes y^{\gamma_2} \wedge z^{\gamma_3} \\
+ \Omega'_2 x^{\alpha_1 - 1} y^{\alpha_2} z^{\alpha_3} \otimes x^{\gamma_1} \wedge z^{\gamma_3},
\end{aligned}$$

avec $\Omega'_1$ et $\Omega'_2$ définis en (12) et (14). L'espace $(K_C)_*$ est donc réduit à 0 pour $* > 3$, et on n'a pas d'homologie en degré supérieur à 3.

Le cas $\lambda$ non racine de l'unité étant un cas particulier du cas libre étudié à la section précédente, on se limite au cas où $\lambda$ est d'ordre $n$ dans $k^*$ pour un entier $n \geq 1$. Remarquons que si $n = 1$ alors $\lambda = 1$ et $A = A_1(k) \otimes k[z]$.

## 5.1 Homologie en degré 3.

Calculons $H_3(K_C) = \operatorname{Ker} d_3 / \operatorname{Im} d_3$. Puisque $(K_C)_4 = 0$, on a $H_3(K_C) = \operatorname{Ker} d_3$. En degré 3 le complexe $K_C$ est :

$$(K_C)_3 = \bigoplus_{\lambda^{\alpha_3+1}=\lambda^{\alpha_2-\alpha_1}=1} k.x^{\alpha_1} y^{\alpha_2} z^{\alpha_3} \otimes x \wedge y \wedge z.$$

La différentielle s'exprime comme suit en degré 3 :

$$\begin{aligned} d_3(x^{\alpha_1} y^{\alpha_2} z^{\alpha_3} \otimes x \wedge y \wedge z) = & -\lambda\alpha_2 x^{\alpha_1} y^{\alpha_2-1} z^{\alpha_3} \otimes y \wedge z \\ & -\lambda\alpha_1 x^{\alpha_1-1} y^{\alpha_2} z^{\alpha_3} \otimes x \wedge z. \end{aligned} \tag{23}$$

**Proposition 5.1.1** *Supposons $\lambda$ d'ordre $n$ dans $k^*$. Alors*

$$H_3(K_C) = \bigoplus_{s \geq 1} k.z^{sn-1} \otimes x \wedge y \wedge z.$$

*C'est donc un module de dimension infinie sur $k$, et de dimension 1 sur le centre de $A$.*

*Preuve.* Puisque les monômes forment une base de PBW, on déduit de la formule (23) que

$$\operatorname{Ker} d_3 = \bigoplus_{\lambda^{\alpha_3+1}=1} k.z^{\alpha_3} \otimes x \wedge y \wedge z,$$

ce qui démontre le résultat. □

## 5.2 Homologie en degré 2.

En degré différentiel 2 le complexe $(K_C)_2$ vaut :

$$\begin{aligned} (K_C)_2 \;=\; & \bigoplus_{\lambda^{\alpha_3+1}=\lambda^{\alpha_2-\alpha_1+1}=1} k.x^{\alpha_1} y^{\alpha_2} z^{\alpha_3} \otimes y \wedge z \\ & \oplus \bigoplus_{\lambda^{\beta_3+1}=\lambda^{\beta_2-\beta_1-1}=1} k.x^{\beta_1} y^{\beta_2} z^{\beta_3} \otimes x \wedge z \\ & \oplus \bigoplus_{\left\{\begin{array}{l} \lambda^{\gamma_3}=1 \\ \gamma_3 = 0 \text{ ou } \lambda^{\gamma_2-\gamma_1}=1 \end{array}\right.} k.x^{\gamma_1} y^{\gamma_2} z^{\gamma_3} \otimes x \wedge y, \end{aligned}$$

et sa différentielle est donnée par :

$$d_2(x^{\alpha_1} y^{\alpha_2} z^{\alpha_3} \otimes y \wedge z) = \alpha_1 \lambda x^{\alpha_1-1} y^{\alpha_2} z^{\alpha_3} \otimes z, \tag{24}$$

$$d_2(x^{\beta_1} y^{\beta_2} z^{\beta_3} \otimes x \wedge z) = -\beta_2 \lambda x^{\beta_1} y^{\beta_2-1} z^{\beta_3} \otimes z, \tag{25}$$

$$d_2(x^{\gamma_1} y^{\gamma_2} z^{\gamma_3} \otimes x \wedge y) = -\gamma_2 x^{\gamma_1} y^{\gamma_2-1} z^{\gamma_3} \otimes y - \gamma_1 x^{\gamma_1-1} y^{\gamma_2} z^{\gamma_3} \otimes x. \tag{26}$$

Soit $X \in \operatorname{Ker} d_2$, on l'écrit

$$X = \sum a_\alpha x^{\alpha_1} y^{\alpha_2} z^{\alpha_3} \otimes y \wedge z + \sum b_\beta x^{\beta_1} y^{\beta_2} z^{\beta_3} \otimes x \wedge z + \sum c_\gamma x^{\gamma_1} y^{\gamma_2} z^{\gamma_3} \otimes x \wedge y,$$

avec $a_\alpha, b_\beta$ et $c_\gamma$ dans $k$ pour tous $\alpha, \beta, \gamma$. Alors $d_2(X) = 0$ implique les deux relations suivantes :

- $c_\gamma = 0$ dès que $\gamma_1 + \gamma_2 \geq 1$ ;
- $-(\alpha_1+1) a_{\alpha_1+1,\alpha_2,\alpha_3} + (\alpha_2+1) b_{\alpha_1,\alpha_2+1,\alpha_3} = 0$ pour tout $\alpha \in \mathbb{N}^3$.

**Remarque :** La formule (23) implique que l'image de $d_3$ est incluse dans la somme directe $A \otimes x \wedge z \;\oplus\; A \otimes y \wedge z$.

**Proposition 5.2.1** *Supposons $\lambda$ d'ordre $n$ dans $k^*$. Alors*

$$H_2(K_C) = \bigoplus_{s \geq 0} k.z^{sn} \otimes x \wedge y.$$

*C'est donc un module de dimension infinie sur $k$, et de dimension 1 sur le centre de $A$.*

*Preuve.* Il suffit de montrer que dans $X$, la partie

$$X_1 = \sum a_\alpha x^{\alpha_1} y^{\alpha_2} z^{\alpha_3} \otimes y \wedge z + \sum b_\beta x^{\beta_1} y^{\beta_2} z^{\beta_3} \otimes x \wedge z$$

possède un antécédent pour $d_2$ dans $(K_C)_3$. Posons :

$$\begin{array}{ll} d_\alpha = -(\lambda\alpha_2)^{-1} a_{\alpha_1,\alpha_2-1,\alpha_3} & \text{pour tous } \alpha_2 \geq 1 \text{ et } \alpha_1 \geq 0; \\ d_\alpha = -(\lambda\alpha_1)^{-1} a_{\alpha_1-1,\alpha_2,\alpha_3} & \text{pour } \alpha_2 = 0 \text{ et } \alpha_1 \geq 1. \end{array}$$

Alors $Y_1 = \sum d_\alpha . x^{\alpha_1} y^{\alpha_2} z^{\alpha_3} \otimes x \wedge y \wedge z$ est un élément de $(K_C)_3$, et on vérifie facilement que $d_3(Y_1) = X_1$. □

## 5.3 Homologie en degré 1.

En degré différentiel 1 le complexe $(K_C)_1$ vaut :

$$\begin{aligned} (K_C)_1 \quad = \quad & \bigoplus_{\alpha \in C_1} k.x^{\alpha_1} y^{\alpha_2} z^{\alpha_3} \otimes z \\ & \oplus \bigoplus_{\beta \in C_2} k.x^{\beta_1} y^{\beta_2} z^{\beta_3} \otimes y \\ & \oplus \bigoplus_{\gamma \in C_3} k.x^{\gamma_1} y^{\gamma_2} z^{\gamma_3} \otimes x, \end{aligned}$$

avec :

- $C_1 = \{\alpha \in \mathbb{N}^3 |\ [\lambda^{\alpha_3+1} = 1 \text{ ou } \alpha_1 = \alpha_2 = 0], \text{ et } \lambda^{\alpha_2-\alpha_1} = 1\}$, c'est-à-dire $C_1 = \{\alpha \in \mathbb{N}^3 |\ [\alpha_3 \equiv -1(n) \text{ ou } \alpha_1 = \alpha_2 = 0], \text{ et } \alpha_2 - \alpha_1 \equiv 0(n)\}$ ;

– $C_2 = \{\beta \in \mathbb{N}^3 |\ \lambda^{\beta_3} = 1,\ \text{et}\ [\lambda^{\beta_2-\beta_1+1} = 1\ \text{ou}\ \beta_3 = 0]\}$, c'est-à-dire $C_2 = \{\beta \in \mathbb{N}^3 |\ \beta_3 \equiv 0(n),\ \text{et}\ [\beta_2 - \beta_1 \equiv -1(n)\ \text{ou}\ \beta_3 = 0]\}$ ;

– $C_3 = \{\gamma \in \mathbb{N}^3 |\ \lambda^{\gamma_3} = 1,\ \text{et}\ [\gamma_3 = 0\ \text{ou}\ \lambda^{\gamma_2-\gamma_1-1} = 1]\}$, c'est-à-dire $C_3 = \{\gamma \in \mathbb{N}^3 |\ \gamma_3 \equiv 0(n),\ \text{et}\ [\gamma_3 = 0\ \text{ou}\ \gamma_2 - \gamma_1 \equiv 1(n)]\}$.

La différentielle en degré 1 vaut :

$$\begin{aligned} d_1(x^{\alpha_1}y^{\alpha_2}z^{\alpha_3} \otimes z) &= 0, \\ d_1(x^{\beta_1}y^{\beta_2}z^{\beta_3} \otimes y) &= \beta_1 x^{\beta_1-1}y^{\beta_2}z^{\beta_3} \otimes 1, \\ d_1(x^{\gamma_1}y^{\gamma_2}z^{\gamma_3} \otimes x) &= -\gamma_2 x^{\gamma_1}y^{\gamma_2-1}z^{\gamma_3} \otimes 1. \end{aligned}$$

Soit $X \in \operatorname{Ker} d_1$, on l'écrit

$$X = \sum a_\alpha x^{\alpha_1}y^{\alpha_2}z^{\alpha_3} \otimes z + \sum b_\beta x^{\beta_1}y^{\beta_2}z^{\beta_3} \otimes y + \sum c_\gamma x^{\gamma_1}y^{\gamma_2}z^{\gamma_3} \otimes x,$$

avec $a_\alpha, b_\beta$ et $c_\gamma$ dans $k$ pour tous $\alpha, \beta, \gamma$. Alors $d_1(X) = 0$ implique :

$$(\alpha_1 + 1)b_{\alpha_1+1,\alpha_2,\alpha_3} - (\alpha_2 + 1)c_{\alpha_1,\alpha_2+1,\alpha_3} = 0 \text{ pour tout } \alpha \in \mathbb{N}^3.$$

**Proposition 5.3.1** *Supposons $\lambda$ d'ordre $n$ dans $k^*$. Alors*

$$H_1(K_C) = \bigoplus_{s \geq 0} \bigoplus_{0 \leq t \leq n-2} k.z^{sn+t} \otimes z.$$

*C'est donc un module de dimension infinie sur $k$ si $\lambda \neq 1$, et de dimension $n - 1$ sur le centre de $A$.*

*Preuve.* 1. Montrons que la partie $X_1 = \sum b_\beta x^{\beta_1}y^{\beta_2}z^{\beta_3} \otimes y + \sum c_\gamma x^{\gamma_1}y^{\gamma_2}z^{\gamma_3} \otimes x$ de $X$ est dans l'image de $d_2$. Comme pour la proposition précédente, on pose $d_\alpha = -\alpha_1^{-1}c_{\alpha_1-1,\alpha_2,\alpha_3}$ pour $\alpha_1 \geq 1$, et $d_\alpha = -\alpha_2^{-1}b_{0,\alpha_2-1,\alpha_3}$ pour $\alpha_1 = 0$ et $\alpha_2 \geq 1$. La formule (26) montre qu'alors $d_2(\sum d_\alpha x^{\alpha_1}y^{\alpha_2}z^{\alpha_3} \otimes x \wedge y) = X_1$.
2. Reste à considérer l'image par $d_2$ des polynômes en $x^{\alpha_1}y^{\alpha_2}z^{\alpha_3} \otimes y \wedge z$ et $x^{\beta_1}y^{\beta_2}z^{\beta_3} \otimes x \wedge z$. Or les formules (24) et (25) donnent :

$$\begin{aligned} d_2\left(\sum e_\alpha x^{\alpha_1}y^{\alpha_2}z^{\alpha_3} \otimes y \wedge z + \sum f_\alpha x^{\alpha_1}y^{\alpha_2}z^{\alpha_3} \otimes x \wedge z\right) = \\ \sum(\alpha_1 e_\alpha \lambda x^{\alpha_1-1}y^{\alpha_2}z^{\alpha_3} \otimes z) + \sum(-\beta_2 f_\beta \lambda x^{\beta_1}y^{\beta_2-1}z^{\beta_3} \otimes z). \end{aligned}$$

Pour apparaître dans l'expression de ce polynôme, un monôme $x^{\alpha_1}y^{\alpha_2}z^{\alpha_3} \otimes z$ dans $K_C$, donc tel que $(\alpha_1, \alpha_2, \alpha_3 + 1) \in C$, doit vérifier $(\alpha_1 + 1, \alpha_2 + 1, \alpha_3 + 1) \in C$. Réciproquement, si $(\alpha_1 + 1, \alpha_2 + 1, \alpha_3 + 1) \in C$, alors

$$d_2(x^{\alpha_1+1}y^{\alpha_2}z^{\alpha_3} \otimes y \wedge z) = (\alpha_1 + 1)x^{\alpha_1}y^{\alpha_2}z^{\alpha_3} \otimes z.$$

Soit $(\alpha_1, \alpha_2, \alpha_3)$ un triplet tel que $(\alpha_1, \alpha_2, \alpha_3 + 1) \in C$. Cherchons la condition nécessaire et suffisante pour avoir $(\alpha_1+1, \alpha_2+1, \alpha_3+1) \notin C$. Puisque $\alpha_3+1 \neq 0$, on a nécessairement $\lambda^{\alpha_1-\alpha_2} = 1$, d'où $\lambda^{\alpha_1+1-(\alpha_2+1)} = 1$. Par ailleurs, soit $\lambda^{\alpha_3+1} = 1$, soit $\alpha_1 = \alpha_2 = 0$. Si $\lambda^{\alpha_3+1} = 1$, alors $(\alpha_1 + 1, \alpha_2 + 1, \alpha_3 + 1) \in C$. La seule possibilité est donc que $\alpha_1 = \alpha_2 = 0$, et $\lambda^{\alpha_3+1} \neq 1$, d'où le résultat. □

## 5.4 Homologie en degré 0.

On doit calculer $H_0(K_C) = \operatorname{Ker} d_0 / \operatorname{Im} d_1 = \operatorname{Coker} d_1$.

**Proposition 5.4.1** *Supposons $\lambda$ d'ordre $n$ dans $k^*$. Alors*

$$H_0(K_C) = \bigoplus_{s \geq 0} \bigoplus_{1 \leq t \leq n-1} k.z^{sn+t} \otimes 1.$$

*C'est donc un module de dimension infinie sur $k$ si $\lambda \neq 1$, et de dimension $n-1$ sur le centre de $A$.*

*Preuve.* Comme dans le cas précédent, une base de $H_0(K_C)$ est constituée des monômes $x^{\alpha_1} y^{\alpha_2} z^{\alpha_3} \otimes 1$ tels que $(\alpha_1, \alpha_2, \alpha_3) \in C$ et $(\alpha_1 + 1, \alpha_2 + 1, \alpha_3) \notin C$. A nouveau ceci n'est possible que si $\alpha_1 = \alpha_2 = 0$ et $\lambda^{\alpha_3} \neq 1$. □

## 5.5 Homologie de Hochschild de l'algèbre mixte minimale.

On peut maintenant synthétiser les résultats de cette section.

**Théorème 5.5.1** *Soit $\lambda$ d'ordre $n$ dans $k^*$, et $A = A^{\lambda}_{2,1}$ la $k$-algèbre engendrée par $x, y, z$ avec les relations*

$$xy - yx = 1, \; xz = \lambda zx, \; yz = \lambda^{-1} zy.$$

*Alors les groupes d'homologie de $A$ sont donnés par :*

- $HH_0(A) = \bigoplus_{s \geq 0} \bigoplus_{1 \leq t \leq n-1} k.z^{sn+t} \otimes 1$, *module de dimension $n-1$ sur le centre $k[z^n]$ de $A$ ;*
- $HH_1(A) = \bigoplus_{s \geq 0} \bigoplus_{0 \leq t \leq n-2} k.z^{sn+t} \otimes z$, *module de dimension $n-1$ sur le centre $k[z^n]$ de $A$ ;*
- $HH_2(A) = \bigoplus_{s \geq 0} k.z^{sn} \otimes x \wedge y$, *module de dimension 1 sur le centre $k[z^n]$ de $A$ ;*
- $HH_3(A) = \bigoplus_{s \geq 1} k.z^{sn-1} \otimes x \wedge y \wedge z$, *module de dimension 1 sur le centre $k[z^n]$ de $A$ ;*
- $HH_*(A) = 0$ *si* $* \geq 4$.

*Preuve.* On a remarqué en début de section l'absence d'homologie en degré supérieur à 3. On applique ensuite les résultats des propositions 5.1.1, 5.2.1, 5.3.1, et 5.4.1. □

**Remarque :** Dans le cas où $\lambda = 1$, c'est-à-dire $n = 1$, l'algèbre $A$ est le produit tensoriel de $A_1(k)$ par $k[z]$, et on retrouve l'homologie de Hochschild donnée par la formule de Künneth.

Donnons à titre de comparaison les résultats pour $\lambda$ non racine de l'unité.

**Théorème 5.5.2** *Soit $\lambda \in k^*$ non racine de l'unité, et $A = A_{2,1}^{\lambda}$ la $k$-algèbre engendrée par $x, y, z$ avec les relations*

$$xy - yx = 1,\ xz = \lambda zx,\ yz = \lambda^{-1} zy.$$

*Alors les groupes d'homologie de $A$ sont donnés par :*

- $HH_0(A) = \bigoplus_{s\geq 1} k.z^s \otimes 1$ ;
- $HH_1(A) = \bigoplus_{s\geq 0} k.z^s \otimes z$ ;
- $HH_2(A) = k.1 \otimes x \wedge y$ ;
- $HH_*(A) = 0$ *si* $* \geq 3$.

*Preuve.* Comme on l'a remarqué en début de section ce cas correspond au cas libre. On applique donc le théorème 4.4.1 avec $r = 1$ et $n = 2$. □

**Remarques :** • Les théorèmes 5.5.1 et 5.5.2 font apparaître que dans le cas mixte minimal l'homologie dépend du paramètre de quantification $\lambda$. Dans le cas semi-classique où $n = r$, le théorème 3.2.4 montre *a contrario* que l'homologie est indépendante de la matrice $\Lambda$ des paramètres.
• L'algèbre $A$ présentée dans cette section est une "Generalized Weyl algebra" au sens de Bavula (voir [4]), mais pas du type étudié dans [8] par Farinati, Solotar et Suarez-Alvarez. En effet dans leur étude ces auteurs considèrent des algèbres GWA construites à partir de l'algèbre de polynômes $k[h]$, avec pour automorphisme de $k[h]$ traduisant les relations de commutation une translation $h \mapsto h + a$ alors que l'algèbre $A$ considérée ici est construite à partir de $k[h, z]$ à l'aide d'un automorphisme $\sigma$ mêlant une translation $h \mapsto h + 1$ et une *homothétie* $z \mapsto \lambda z$.

# 6 Généralité sur la dualité pour les algèbres $A_{n,r}^{\Lambda}$.

Afin de démontrer dans la prochaine section un résultat de dualité pour les algèbres semi-classiques, on va dans la section présente établir les liens existant entre l'homologie et la cohomologie de Hochschild pour une algèbre $A_{n,r}^{\Lambda}$. Nous allons pour cela décrire le complexe calculant la cohomologie de $A_{n,r}^{\Lambda}$ découlant de la résolution présentée à la proposition 1.5.2, et expliciter un certain nombre d'isomorphismes de complexes.

## 6.1 Cohomologie de Hochschild des algèbres $A_{n,r}^{\Lambda}$.

Rappelons les notations. On se donne une matrice $\Lambda$ multiplicativement antisymétrique de taille $n$, et on définit $A_{n,r}^{\Lambda}$ conformément à 1.1. Par la proposition 1.7.1, $A_{n,r}^{\Lambda}$ est une algèbre enveloppante quantique $U_Q(V, f)$, où $V$ est un $k$-espace vectoriel de base $(x_1, \ldots, x_r, y_1, \ldots, y_n)$, où $Q = Q(\Lambda)$ est la matrice $Q(\Lambda) = (\widetilde{\lambda}_{i,j})_{1\leq i,j\leq n+r} = \begin{pmatrix} \Lambda_r & \Lambda_{r,n}^{-1} \\ \Lambda_{n,r}^{-1} & \Lambda \end{pmatrix}$ introduite au paragraphe 2.1, et $f$ est la forme linéaire de $V \otimes V$ valant 1 en $x_i \otimes y_i$, valant $-1$ en $y_i \otimes x_i$ et 0 sur les autres éléments de base. On note $U = A^{\Lambda}_{n,r} = U_Q(V, f)$.
L'espace $\Lambda_Q V$ défini au paragraphe 1.5 admet pour base la famille $(x_1^{\gamma_1} \wedge \ldots y_n^{\delta_n})$. Renommons $v_1, \ldots, v_{n+r}$ les générateurs de $V$, en posant $v_i = x_i$ pour $i \leq r$ et $v_i = y_{i-r}$ pour $i \geq r+1$. Alors pour tout $* \in \mathbb{N}$ on a

$$\Lambda^*_Q V = \bigoplus_{i_1 < \ldots < i_*} k.v_{i_1} \wedge \ldots \wedge v_{i_*}.$$

Enfin, si $\{i_1, \ldots, i_*\}$ est une famille d'indices tels que $1 \leq i_1 < \ldots < i_* \leq r+n$, on note $\{j_t\}_{1 \leq t \leq n+r-*} = \overline{\{i_s\}}_{1 \leq s \leq *}$ l'ensemble des indices complémentaires tels que $1 \leq j_1 < \ldots < j_{n+r-*} \leq n+r$ et $\{i_s\} \cup \{j_t\} = \{1, \ldots, n+r\}$.
Puisque $U \otimes \Lambda^*_Q V \otimes U$ est une résolution libre de $U$ par des $U^e$-bimodules, on en déduit un complexe $(K^*, {}^t\partial)$ calculant la cohomologie de Hochschild de $U$. Comme $k$-espace vectoriel,

$$K^* = \mathrm{Hom}_{U^e}(U \otimes \Lambda^*_Q V \otimes U, U),$$

et la différentielle est la transposée de la différentielle

$$\partial : U \otimes \Lambda^{*+1}_Q V \otimes U \to U \otimes \Lambda^*_Q V \otimes U,$$

à savoir :

$$\begin{array}{rl} {}^t\partial : & K^* \to K^{*+1} \\ & \varphi \mapsto \varphi \circ d. \end{array}$$

**Proposition 6.1.1** *Soient $r \leq n$ deux entiers, $\Lambda$ une matrice multiplicative antisymétrique de taille $n$, et $A^{\Lambda}_{n,r}$ l'algèbre définie en 1.1. Alors :*

$$HH^*(A^{\Lambda}_{n,r}) = H^*(K, {}^t\partial).$$

*Preuve.* C'est la définition de la cohomologie calculée à partir de la résolution présentée à la proposition 1.5.2. □

**Corollaire 6.1.2** *Les modules de cohomologie de Hochschild de $A^{\Lambda}_{n,r}$ sont tous nuls en degré strictement supérieur à $n+r$.*

*Preuve.* Les espaces vectoriels $K^*$ sont tous nuls pour $* > n+r$. □

Il existe entre $\mathrm{Hom}_{U^e}(U \otimes \Lambda^*_Q V \otimes U, U)$ et $U \otimes \Lambda^{n+r-*}_Q V$ une suite d'isomorphismes naturels d'espaces vectoriels, qui envoient successivement $\mathrm{Hom}_{U^e}(U \otimes \Lambda^*_Q V \otimes U, U)$ sur $\mathrm{Hom}_k(\Lambda^*_Q V, U)$, puis $\mathrm{Hom}_k(\Lambda^*_Q V, U)$ sur $U \otimes (\Lambda^*_Q V)'$, et enfin $U \otimes (\Lambda^*_Q V)'$ sur $U \otimes \Lambda^{n+r-*}_Q V$. Afin de savoir sous quelle condition on peut en déduire une dualité entre l'homologie et la cohomologie de Hochschild de $U$, nous allons étudier sous quelles conditions ces isomorphismes d'espace vectoriel peuvent être des isomorphismes de complexes différentiels.

## 6.2 Conjugaison de structures différentielles.

Suivant le principe qui a permis de passer de la proposition 1.5.2 à la proposition 1.5.3, explicitons un résultat général, qui va nous permettre de transporter les structures différentielles de $\mathrm{Hom}_{U^e}(U \otimes \Lambda_Q^* V \otimes U, U)$ et $U \otimes \Lambda^{n+r-*}V$ respectivement sur $\mathrm{Hom}_k(\Lambda_Q^* V, U)$ et $U \otimes (\Lambda_Q^* V)'$, afin de les comparer.

**Lemme 6.2.1 (Lemme de conjugaison)** *Soient $(C_*, d)$ un $k$-complexe différentiel, et $M_*$ un $k$-espace vectoriel gradué, tel qu'il existe un isomorphisme $\Phi$ d'espace vectoriels gradués de degré 0, de source $C_*$ et de but $M_*$. Alors l'application $\widetilde{d} = \Phi \circ d \circ \Phi^{-1}$ est telle que $\widetilde{d}^2 = 0$, et $(M_*, \widetilde{d})$ est un complexe différentiel. L'application $\Phi$ est alors un isomorphisme de complexes, et on a :*

$$H_*(C, d) = H_*(M, \widetilde{d}).$$

*Preuve.* Il suffit de vérifier que par construction $\Phi$ est un isomorphisme de complexes. □

Appliquons ce résultat au complexe $(K^*, {}^t\partial)$ décrit ci-dessus, et au complexe $(K_*, d)$ décrit à la fin du paragraphe 1.7 par les formules (9) à (14).
Il existe un $k$-isomorphisme $\Phi_{1,*}$ de $\mathrm{Hom}_{U^e}(U \otimes \Lambda^* V \otimes U, U)$ sur $\mathrm{Hom}_k(\Lambda_Q^* V, U)$, défini pour $\varphi \in \mathrm{Hom}_{U^e}(U \otimes \Lambda^* V \otimes U, U)$ par

$$\Phi_{1,*}(\varphi)(v_{i_1} \wedge \ldots \wedge v_{i_*}) = \varphi(1 \otimes v_{i_1} \wedge \ldots \wedge v_{i_*} \otimes 1).$$

On calcule alors le conjugué $D = \Phi_{1,*+1} \circ {}^t\partial \circ \Phi_{1,*}^{-1}$ de ${}^t\partial$ par $\Phi_1$. Posons

$$L^* = \mathrm{Hom}_k(\Lambda_Q^* V, U),$$

alors on a

$$\begin{aligned} D:\ & L^* && \to L^{*+1} \\ & \varphi && \mapsto D(\varphi), \end{aligned}$$

où $D(\varphi)$ est définie par :

$$\begin{aligned} D(\varphi)(v_{i_1} \wedge \ldots \wedge v_{i_{*+1}}) = \\ & \sum_{k=1}^{*+1} (-1)^{k-1} \Big( (\prod_{s<k} \widetilde{\lambda}_{i_s,i_k}) v_{i_k} \varphi(v_{i_1} \wedge \ldots \widehat{v}_{i_k} \ldots \wedge v_{i_{*+1}}) \\ & \qquad - (\prod_{s>k} \widetilde{\lambda}_{i_k,i_s}) \varphi(v_{i_1} \wedge \ldots \widehat{v}_{i_k} \ldots \wedge v_{i_{*+1}}) v_{i_k} \Big). \end{aligned} \tag{27}$$

**Proposition 6.2.2** *Le complexe $(K^*, {}^t\partial)$ de la proposition 6.1.1 est isomorphe au complexe $(L^*, D)$ décrit ci-dessus.*

*Preuve.* En vertu du lemme 6.2.1, ceci découle de la définition de $(L^*, D)$. □

On a donc un diagramme commutatif :

$$\begin{array}{ccc} K^* & \xrightarrow{\quad {}^t\partial \quad} & K^{*+1} \\ \Phi_{1,*}\Big\downarrow & \circlearrowleft & \Big\downarrow \Phi_{1,*+1} \\ \mathrm{Hom}(\Lambda_Q^* V, U) & \xrightarrow{\quad D \quad} & \mathrm{Hom}(\Lambda_Q^{*+1} V, U) \end{array}$$

Pour calculer la cohomologie de Hochschild d'une algèbre $A_{n,r}^{\Lambda}$ on utilise à présent le complexe $(L^*, D_1)$. C'est notamment avec ce complexe qu'on calculera la cohomologie de l'algèbre mixte minimale à la section 8.

Etablissons maintenant un résultat de même nature pour le complexe $(K_*, d)$ défini en 1.7 calculant l'homologie de Hochschild de l'algèbre $U$. On rappelle qu'en tant qu'espace vectoriel, $K_* = U \otimes \Lambda_Q^* V$, et la différentielle $d$ est donnée par la formule (10) à la page 11.

**Lemme 6.2.3** *L'application canonique* $\psi_* : \Lambda_Q^* V \otimes \Lambda_Q^{n+r-*} V \to k \otimes v_1 \wedge \ldots \wedge v_{n+r}$ *définie par* $\psi_*(v_{i_1} \wedge \ldots \wedge v_{i_*} \otimes v_{j_1} \wedge \ldots \wedge v_{j_{n+r-*}}) = v_{i_1} \wedge \ldots \wedge v_{i_*} \wedge v_{j_1} \wedge \ldots \wedge v_{j_{n+r-*}}$, *induit un isomorphisme* $\overline{\psi}_* : \Lambda_Q^{n+r-*} V \to (\Lambda_Q^* V)'$ *défini par*

$$\overline{\psi}_*(v_{j_1} \wedge \ldots \wedge v_{j_{n+r-*}}) = \psi(. \otimes v_{j_1} \wedge \ldots \wedge v_{j_{n+r-*}}).$$

*Preuve.* En fait, $\overline{\psi}_*$ n'est autre que l'application linéaire envoyant l'élément de base $v_{j_1} \wedge \ldots \wedge v_{j_{n+r-*}}$ sur $\Theta_*(i_1, \ldots, i_*)(v_{i_1} \wedge \ldots \wedge v_{i_*})'$, où $\{i_1, \ldots, i_*\}$ est le $*$-uplet complémentaire de $\{j_1, \ldots, j_{n+r-*}\}$, et $\Theta_*(i_1, \ldots, i_*)$ est l'élément de $k^*$ défini par :

$$\Theta_*(i_1, \ldots, i_*) = \prod_{k<i_*,\ k\notin\{i_s\}} (-\widetilde{\lambda}_{i_*,k}) \prod_{k<i_{*-1},\ k\notin\{i_s\}} (-\widetilde{\lambda}_{i_{*-1},k}) \ldots \prod_{k<i_1|k\notin\{i_s\}} (-\widetilde{\lambda}_{i_1,k}). \tag{28}$$

Cette application linéaire envoie donc une base de $\Lambda_Q^{n+r-*} V$ sur une base de $(\Lambda_Q^* V)'$ : la base duale de la base de $\Lambda_Q^* V$ constituée des $v_{i_1} \wedge \ldots \wedge v_{i_*}$. $\square$

L'isomorphisme $\overline{\psi}_*$ induit un autre isomorphisme $\Phi_{2,*} = id \otimes \overline{\psi}_*$ de $U \otimes \Lambda_Q^{n+r-*} V$ sur $U \otimes (\Lambda_Q^* V)'$. Or $U \otimes \Lambda_Q^{n+r-*} V = K_{n+r-*}$, et on va donc définir une différentielle $\Delta$ sur le complexe $U \otimes (\Lambda_Q V)'$, telle que le diagramme suivant commute :

$$\begin{array}{ccc} U \otimes (\Lambda_Q^* V)' & \xrightarrow{\quad \Delta \quad} & U \otimes (\Lambda_Q^{*+1} V)' \\ \Phi_{2,*}\Big\uparrow & \circlearrowleft & \Big\uparrow \Phi_{2,*+1} \\ K_{n+r-*} & \xrightarrow{\quad d \quad} & K_{n+r-*-1} \end{array}$$

La différentielle $\Delta$, qui n'est autre que $\Phi_{2,*+1} \circ d \circ \Phi_{2,*}^{-1}$, se calcule de la façon suivante : soit $(i_1, \ldots, i_*)$, un $*$-uplet, avec $i_1 < \ldots < i_*$, on note $j_1, \ldots, j_{n+r-*}$ les éléments complémentaires : $\{j_1, \ldots, j_{n+r-*}\} = \overline{\{i_1, \ldots, i_*\}}$. Alors

$$\begin{aligned} \Delta(a \otimes (v_{i_1} \wedge \ldots \wedge v_{i_*})') = {} & \Theta_*^{-1}(i_1, \ldots, i_*) \times \\ & \sum_{k=1}^{n+r-*} (-1)^{k-1} \Big( (\prod_{s<k} \widetilde{\lambda}_{j_s,j_k}) a v_{j_k} - (\prod_{s>k} \widetilde{\lambda}_{j_k,j_s}) v_{j_k} a \Big) \\ & \otimes \Theta_{*+1}(i_1, \ldots, j_k, \ldots, i_*)(v_{i_1} \wedge \ldots v_{j_k} \ldots \wedge v_{i_*})'. \end{aligned} \tag{29}$$

**Proposition 6.2.4** *Le complexe $(U \otimes (\Lambda_Q^* V)', \Delta)$ défini ci-dessus a pour homologie :*

$$H^*(U \otimes (\Lambda_Q V)') = HH_{n+r-*}(U).$$

*Preuve.* Par définition de $\Delta$ et par le lemme 6.2.1 on a pour tout $* \in \mathbb{N}$ l'égalité $H^*(U \otimes (\Lambda_Q^* V)', \Delta) = H_{n+r-*}(K_*, d)$. On conclut à l'aide de la proposition 1.7.3. $\square$

## 6.3 Lien entre l'homologie et la cohomologie des algèbres $A_{n,r}^{\Lambda}$.

On a transporté les structures de complexes différentiels de $(K_*, d)$ et $(K^*, {}^t\partial)$ sur les deux espaces $\mathrm{Hom}(\Lambda_Q^* V, U)$ et $U \otimes (\Lambda_Q^* V)'$. A nouveau ces deux espaces sont liés par un isomorphisme naturel. On va comparer leurs structures de complexes différentiels à l'aide de cet isomorphisme.
Appelons $\Phi_{3,*}$ l'isomorphisme de $U \otimes (\Lambda_Q^* V)'$ sur $\mathrm{Hom}(\Lambda_Q^* V, U)$, défini par :

$$\Phi_{3,*}(a \otimes \varphi)(v_{i_1} \wedge \ldots \wedge v_{i_*}) = \varphi(v_{i_1} \wedge \ldots \wedge v_{i_*})a. \qquad (30)$$

Considérons alors le diagramme suivant :

$$\begin{array}{ccc}
\mathrm{Hom}(\Lambda_Q^* V, U) & \xrightarrow{D} & \mathrm{Hom}(\Lambda_Q^{*+1} V, U) \\
\Phi_{3,*} \uparrow & & \uparrow \Phi_{3,*+1} \\
U \otimes (\Lambda_Q^* V)' & \xrightarrow{\Delta} & U \otimes (\Lambda_Q^{*+1} V)'
\end{array}$$

Ce diagramme n'est pas commutatif *a priori*. On va donc préciser les valeurs de $\Phi_{3,*+1} \circ \Delta$ et $D \circ \Phi_{3,*}$.

**Lemme 6.3.1** *Soient $a \in U$, et $1 \leq i_1 < \ldots < i_* \leq n + r$. Soit alors le complémentaire $\{j_1, \ldots, j_{n+r-*}\} = \overline{\{i_s\}}$. Alors $\Phi_{3,*+1} \circ \Delta(a \otimes (v_{i_1} \wedge \ldots \wedge v_{i_*})')$ est l'application linéaire qui envoie un élément $v_{\alpha_1} \wedge \ldots \wedge v_{\alpha_{*+1}}$ de $\Lambda_Q^{*+1} V$ sur :*

$$\sum_{k=1}^{n+r-*} (-1)^{k-1} \omega_1(\alpha_1, \ldots, \alpha_{*+1}; k) a v_{j_k} - \omega_2(\alpha_1, \ldots, \alpha_{*+1}; k) v_{j_k} a,$$

*avec*

$$\begin{aligned}
\omega_1(\alpha_1, \ldots, \alpha_{*+1}; k) = \;& \Theta_*^{-1}(i_1, \ldots, i_*) \Theta_{*+1}(i_1, \ldots, j_k, \ldots, i_*) \prod_{s<k} \widetilde{\lambda}_{j_s, j_k} \\
& \textit{si } (\alpha_1, \ldots, \alpha_{*+1}) = (i_1, \ldots, j_k, \ldots, i_*), \\
\omega_1(\alpha_1, \ldots, \alpha_{*+1}; k) = \;& 0 \;\; \textit{sinon};
\end{aligned} \qquad (31)$$

*et*

$$\begin{aligned}
\omega_2(\alpha_1, \ldots, \alpha_{*+1}; k) = \;& \Theta_*^{-1}(i_1, \ldots, i_*) \Theta_{*+1}(i_1, \ldots, j_k, \ldots, i_*) \prod_{s>k} \widetilde{\lambda}_{j_k, j_s} \\
& \textit{si } (\alpha_1, \ldots, \alpha_{*+1}) = (i_1, \ldots, j_k, \ldots, i_*), \\
\omega_2(\alpha_1, \ldots, \alpha_{*+1}; k) = \;& 0 \;\; \textit{sinon}.
\end{aligned} \qquad (32)$$

*Preuve.* Les formules (29) et (30) permettent de calculer l'image d'un élément $v_{\alpha_1} \wedge \ldots \wedge v_{\alpha_{*+1}}$ de $\Lambda_Q^{*+1}V$ par $\Phi_{3,*+1} \circ \Delta(a \otimes (v_{i_1} \wedge \ldots \wedge v_{i_*})')$. Cette image vaut :

$$\Theta_*^{-1}(i_1, \ldots, i_*) \sum_{k=1}^{n+r-*} (-1)^{k-1} \lambda_{*+1}(i_1, \ldots, j_k, \ldots, i_*) \times$$
$$(v_{i_1} \wedge \ldots v_{j_k} \ldots \wedge v_{i_*})'(v_{\alpha_1} \wedge \ldots \wedge v_{\alpha_{*+1}})(\prod_{s<k} \widetilde{\lambda}_{j_s,j_k} a v_{j_k} - \prod_{s>k} \widetilde{\lambda}_{j_k,j_s} v_{j_k} a),$$

ce qui conduit au résultat annoncé. □

**Lemme 6.3.2** *Soient $a \in U$, et $1 \leq i_1 < \ldots < i_* \leq n + r$. Soit alors l'ensemble complémentaire $\{j_1, \ldots, j_{n+r-*}\} = \overline{\{i_s\}}$. Alors $D \circ \Phi_{3,*}(a \otimes (v_{i_1} \wedge \ldots \wedge v_{i_*})')$ est l'application linéaire qui envoie un élément $v_{\alpha_1} \wedge \ldots \wedge v_{\alpha_{*+1}}$ de $\Lambda_Q^{*+1}V$ sur :*

$$\sum_{k=1}^{n+r-*} (-1)^{k-1} \omega_1'(\alpha_1, \ldots, \alpha_{*+1}; k) a v_{j_k} - \omega_2'(\alpha_1, \ldots, \alpha_{*+1}; k) v_{j_k} a, \tag{33}$$

*avec*

$$\omega_1' = (-1)^{*+1} \omega_1,$$
$$\omega_2' = (-1)^{*+1} (\prod_{t=1}^{n+r} \widetilde{\lambda}_{j_k,t}) \omega_2.$$

*Preuve.* Les formules (27) et (30) permettent de calculer l'image d'un élément $v_{\alpha_1} \wedge \ldots \wedge v_{\alpha_{*+1}} \in \Lambda^{*+1}V$ par $D \circ \Phi_{3,*}(a \otimes (v_{i_1} \wedge \ldots \wedge v_{i_*})')$, et on trouve :

$$\begin{aligned}\sum_{t=1}^{*+1} (-1)^{t-1} \Big( & (\prod_{s<t} \widetilde{\lambda}_{\alpha_s,\alpha_t}) v_{\alpha_t} (v_{i_1} \wedge \ldots \wedge v_{i_*})'(v_{\alpha_1} \wedge \ldots \widehat{v}_{\alpha_t} \ldots \wedge v_{\alpha_{*+1}}) a \\ & - (\prod_{s>t} \widetilde{\lambda}_{\alpha_t,\alpha_s})(v_{i_1} \wedge \ldots \wedge v_{i_*})'(v_{\alpha_1} \wedge \ldots \widehat{v}_{\alpha_t} \ldots \wedge v_{\alpha_{*+1}}) a v_{\alpha_t} \Big).\end{aligned} \tag{34}$$

Remarquons tout d'abord que $(v_{i_1} \wedge \ldots \wedge v_{i_*})'(v_{\alpha_1} \wedge \ldots \widehat{v}_{\alpha_t} \ldots \wedge v_{\alpha_{*+1}})$ est non nul si et seulement s'il existe $k$ tel que $\alpha_t = j_k$, et $(\alpha_1, \ldots, \widehat{\alpha}_t, \ldots, \alpha_{*+1}) = (i_1, \ldots, i_*)$. Dans ce cas, $t - 1$ est égal au nombre d'indices $i_s$ strictement inférieurs à $j_k$, nombre que l'on note $|i_s < j_k|$. La formule (34) apparaît alors sous la forme de (33), avec pour coefficients :

$$\begin{aligned}\omega_1'(\alpha_1, \ldots, \alpha_{*+1}; k) = \quad & (-1)^k (-1)^{|i_s<j_k|} \prod_{i_s>j_k} \widetilde{\lambda}_{j_k,i_s} \\ & \quad \text{si } (\alpha_1, \ldots, \alpha_{*+1}) = (i_1, \ldots, j_k, \ldots, i_*), \\ \omega_1'(\alpha_1, \ldots, \alpha_{*+1}; k) = \quad & 0 \ \text{ sinon;}\end{aligned}$$

et

$$\begin{aligned}\omega_2'(\alpha_1, \ldots, \alpha_{*+1}; k) = \quad & (-1)^k (-1)^{|i_s<j_k|} \prod_{i_s<j_k} \widetilde{\lambda}_{i_s,j_k}, \\ & \quad \text{si } (\alpha_1, \ldots, \alpha_{*+1}) = (i_1, \ldots, j_k, \ldots, i_*), \\ \omega_2'(\alpha_1, \ldots, \alpha_{*+1}; k) = \quad & 0 \ \text{ sinon.}\end{aligned}$$

On simplifie l'expression (31) de $\omega_1$ en utilisant le fait que :

$$\Theta_{*+1}(i_1,\ldots,j_k,\ldots,i_*) = \Theta_*(i_1,\ldots,i_*) \prod_{i_s>j_k} (-\widetilde{\lambda}_{i_s,j_k})^{-1} \prod_{t<j_k,\ t\notin\{i_s\}} (-\widetilde{\lambda}_{j_k,t}),$$

ce qui découle de la formule (28). On obtient alors : $\omega_1 = \omega_1' \times (-1)^{*+1}$. De même on simplifie l'expression (32) de $\omega_2$, et on obtient $\omega_2' = (-1)^{*+1}(\prod_{t=1}^{n+r} \widetilde{\lambda}_{j_k,t})\omega_2$. □

# 7 Dualité dans le cas semi-classique.

Dans le cas semi-classique, le conjugué de $\Delta$ par $\Phi_3$ est égal à $D$ à un signe près, ce qui permet d'établir la dualité. Reprenons les notations de la section 3. On se donne un entier $n$, et une matrice multiplicativement antisymétrique $\Lambda$ de taille $n$. Grâce aux résultats de la section précédente on va comparer les modules d'homologie et de cohomologie de $A^{\Lambda}_{n,n}$ en degré $*$ et $2n-*$.

## 7.1 Application de la section 6 au cas semi-classique.

**Théorème 7.1** *Soient $n$ un entier, et $\Lambda \in M_n(k)$ une matrice multiplicativement antisymétrique. Alors l'homologie et la cohomologie de Hochschild de l'algèbre $A^{\Lambda}_{n,n}$ définie en 1.1 avec $r = n$ vérifient :*

$$HH_*(A^{\Lambda}_{n,n}) = HH^{2n-*}(A^{\Lambda}_{n,n}).$$

*Preuve.* Par les propositions 6.2.2 et 6.2.4, on a

$$HH_{2n-*}(A^{\Lambda}_{n,n}) = H^*(U \otimes (\Lambda_Q V)', \Delta) \ \text{ et } \ HH^*(A^{\Lambda}_{n,n}) = H^*(L, D).$$

Dans le cas semi-classique ($n = r$), la formule (7) montre que la $i^{\text{ème}}$ ligne de $Q(\Lambda) = (\widetilde{\lambda}_{i,j})$ contient exactement une fois $\lambda_{i,t}$ et une fois $\lambda^{-1}_{i,t}$ pour tout $t$ dans $\{1,\ldots,n\}$. Dans ce cas, le lemme 6.3.2 implique $D_* = (-1)^{*+1}\Phi_{3,*+1} \circ \Delta_* \circ \Phi^{-1}_{3,*}$. La conjugaison par $\Phi_3$ induit donc un isomorphisme entre les espaces vectoriels $H^*(U \otimes (\Lambda^*_Q V)', \Delta) = \operatorname{Ker}\Delta_*/\operatorname{Im}\Delta_{*-1}$ et $H^*(L, D) = \operatorname{Ker} D_*/\operatorname{Im} D_{*-1}$, ce qui termine la preuve. □

On parle de dualité en homologie pour signifier des résultats du type : "l'algèbre $A$ vérifie $HH_*(A) = HH^{d-*}(A)$ pour un certain entier $d$". Interprétons les résultats précédents dans ce cadre.

## 7.2 Dualité et dimension globale.

**7.2.1** De nombreux cas classiques font apparaître la dimension globale comme indice de dualité. On a vu en 1.3.3 que $A^{\Lambda}_{n,n}$ est un cas particulier d'algèbre $A^{\bar{q},\Lambda}_n(k)$ où $\bar{q} = (q_1,\ldots,q_n)$ ne contient que des 1. On peut donc lui appliquer le théorème 3.8 de [9], d'où il découle que la dimension globale de l'algèbre $A^{\Lambda}_{n,n}$ est $n$. L'algèbre $(A^\Lambda_{n,n})^e = A^\Lambda_{n,n} \otimes (A^\Lambda_{n,n})^{op}$ est encore une algèbre semi-classique, à savoir $A^{\widetilde{\Lambda}}_{2n,2n}$ avec la matrice diagonale par blocs $\widetilde{\Lambda} = \mathrm{Diag}(\Lambda, {}^t\Lambda)$. Elle a donc pour dimension globale $2n$, et si on note $d$ cette dimension globale, le théorème 7.1 s'écrit : $HH_*(A^\Lambda_{n,n}) = HH^{d-*}(A^\Lambda_{n,n})$.

**7.2.2** Pour $n > r$, il n'y a pas d'espoir d'avoir un résultat semblable pour une algèbre $A^\Lambda_{n,r}$ quelconque. En effet dans le cas libre et pour $n > r$, le complexe de Wambst fournit une résolution libre de $(A^\Lambda_{n,r})^e$ de longueur $n+r$. Donc la dimension globale $d$ de $(A^\Lambda_{n,r})^e$ vaut au moins $n + r$. Par ailleurs, lorsque les coefficients $\lambda_{i,j}$ sont indépendants on peut montrer que $\mathcal{Z}(A^\Lambda_{n,r}) = k$ (voir [25]), c'est-à-dire qu'on a $HH^0(A^\Lambda_{n,r}) = k$, tandis que pour tout $d \geq n+r > 2r$ on a vu au théorème 4.4.1 que $HH_d(A^\Lambda_{n,r}) = 0$. On n'a donc pas dans ce cas de dualité avec pour indice la dimension globale de $(A^\Lambda_{n,r})^e$.

**7.2.3** Nous proposerons à la proposition 8.2.1 un exemple pour lequel la dualité n'est jamais vérifiée, même pour un indice autre que la dimension globale de $(A^\Lambda_{n,r})^e$.

# 8 Cohomologie de l'algèbre mixte minimale : un exemple de non dualité.

On reprend les notations de la section 5. On note $V$ le $k$-espace vectoriel de base $(x, y, z)$. Dans tout cette section on note $A = A^\Lambda_{2,1}$ l'algèbre engendré par $x, y, z$ avec les relations

$$xy - yx = 1,\ xz = \lambda zx,\ yz = \lambda^{-1} zy. \tag{35}$$

Les matrices associées à $A$ sont $\Lambda = \begin{pmatrix} 1 & \lambda^{-1} \\ \lambda & 1 \end{pmatrix}$, et $Q = Q(\Lambda) = \begin{pmatrix} 1 & 1 & \lambda \\ 1 & 1 & \lambda^{-1} \\ \lambda^{-1} & \lambda & 1 \end{pmatrix}$.

Si $\lambda = 1$ la formule de Künneth et la dualité de Poincaré pour $A_1$ et $k[z]$ montrent la dualité pour l'algèbre $A$. On suppose donc désormais $\lambda \neq 1$.

## 8.1 Calcul de $HH^0(A)$ et $HH^1(A)$.

On déduit des relations entre les générateurs de $A$ la propriété suivante :

**Lemme 8.1.1** *Dans $A$ pour tout triplet $(\alpha_1, \alpha_2, \alpha_3)$ de $\mathbb{N}^3$ on a les relations suivantes :*

- $[x, x^{\alpha_1} y^{\alpha_2} z^{\alpha_3}] = (1 - \lambda^{-\alpha_3}) x^{\alpha_1+1} y^{\alpha_2} z^{\alpha_3} + \lambda^{-\alpha_3} \alpha_2 x^{\alpha_1} y^{\alpha_2-1} z^{\alpha_3}$ ;
- $[y, x^{\alpha_1} y^{\alpha_2} z^{\alpha_3}] = (1 - \lambda^{\alpha_3}) x^{\alpha_1} y^{\alpha_2+1} z^{\alpha_3} - \alpha_1 x^{\alpha_1-1} y^{\alpha_2} z^{\alpha_3}$ ;
- $[z, x^{\alpha_1} y^{\alpha_2} z^{\alpha_3}] = (\lambda^{-\alpha_1+\alpha_2} - 1) x^{\alpha_1} y^{\alpha_2} z^{\alpha_3+1}$.

*Preuve.* Les calculs se font simplement à partir des relations (35), et sont laissés au lecteur. □

On a alors en degré différentiel 0 et 1 :

$$K^0 = \mathrm{Hom}(k, A),$$

$$K^1 = \mathrm{Hom}(V, A).$$

Pour $\varphi \in K^0$, l'application $D_0(\varphi)$ est définie par :

$$\begin{aligned} D_0(\varphi)(x) &= \; x\varphi(1) - \varphi(1)x, \\ D_0(\varphi)(y) &= \; y\varphi(1) - \varphi(1)y, \\ D_0(\varphi)(z) &= \; z\varphi(1) - \varphi(1)z. \end{aligned}$$

**Lemme 8.1.2** *Le module en cohomologie de degré 0 de $A$ vaut :*

– $HH^0(A) = k$ *si* $\lambda$ *n'est pas racine de l'unité ;*

– $HH^0(A) = k[z^n]$ *si* $\lambda$ *est d'ordre* $n$ *dans* $k^*$.

*Preuve.* L'espace $HH^0(A)$ est isomorphe au centre de $A$, que nous avons déjà décrit au début de la section 5. □

Calculons maintenant le module en cohomologie de degré 1. Rappelons que :

$$K^1 = \mathrm{Hom}(V, A),$$

$$K^2 = \mathrm{Hom}(\Lambda^2_Q V, A).$$

Pour $\varphi \in K^1$, on a :

$$\begin{aligned} D_1(\varphi)(x \wedge y) &= \; x\varphi(y) - \varphi(y)x - (y\varphi(x) - \varphi(x)y), \\ D_1(\varphi)(x \wedge z) &= \; x\varphi(z) - \lambda\varphi(z)x - (\lambda z\varphi(x) - \varphi(x)z), \\ D_1(\varphi)(y \wedge z) &= \; y\varphi(z) - \lambda^{-1}\varphi(z)y - (\lambda^{-1}z\varphi(y) - \varphi(y)z). \end{aligned}$$

**Lemme 8.1.3** *Le module en cohomologie de degré 1 de $A$ vaut :*

– $HH^1(A) = k.z$ *si* $\lambda$ *n'est pas racine de l'unité ;*

– $HH^1(A) = z.k[z^n]$ *si* $\lambda$ *est d'ordre* $n$ *dans* $k^*$.

*Preuve.* On doit déterminer $HH^1(A) = \mathrm{Ker}\, D_1 / \mathrm{Im}\, D_0$. Fixons $\varphi$ dans $\mathrm{Ker}\, D_1$. Notons

$$\varphi(x) = \sum a_\alpha x^{\alpha_1} y^{\alpha_2} z^{\alpha_3}, \;\; \varphi(y) = \sum b_\alpha x^{\alpha_1} y^{\alpha_2} z^{\alpha_3}, \;\; \text{et } \varphi(z) = \sum c_\alpha x^{\alpha_1} y^{\alpha_2} z^{\alpha_3},$$

avec $a_\alpha, b_\beta$ et $c_\gamma$ dans $k$ pour tous $\alpha, \beta, \gamma$. Alors $\varphi$ est un élément de $\mathrm{Ker}\, D_1$ si et seulement si l'on a dans $A$ les relations suivantes :

$$\begin{aligned} &\sum b_\alpha (1 - \lambda^{-\alpha_3}) x^{\alpha_1+1} y^{\alpha_2} z^{\alpha_3} + \sum b_\alpha \lambda^{-\alpha_3} \alpha_2 x^{\alpha_1} y^{\alpha_2-1} z^{\alpha_3} = \\ &\qquad \sum a_\alpha (1 - \lambda^{-\alpha_3}) x^{\alpha_1} y^{\alpha_2+1} z^{\alpha_3} - \sum a_\alpha \alpha_1 x^{\alpha_1-1} y^{\alpha_2} z^{\alpha_3}; \end{aligned} \tag{36}$$

$$\begin{aligned} &\sum c_\alpha (1 - \lambda^{1-\alpha_3}) x^{\alpha_1+1} y^{\alpha_2} z^{\alpha_3} + \sum c_\alpha \lambda^{1-\alpha_3} \alpha_2 x^{\alpha_1} y^{\alpha_2-1} z^{\alpha_3} = \\ &\qquad \sum a_\alpha (\lambda^{1-\alpha_1+\alpha_2} - 1) x^{\alpha_1} y^{\alpha_2} z^{\alpha_3+1}; \end{aligned} \tag{37}$$

$$\begin{aligned} &\sum c_\alpha (1 - \lambda^{-1+\alpha_3}) x^{\alpha_1} y^{\alpha_2+1} z^{\alpha_3} - \sum c_\alpha \alpha_1 x^{\alpha_1-1} y^{\alpha_2} z^{\alpha_3} = \\ &\qquad \sum b_\alpha (\lambda^{-1-\alpha_1+\alpha_2} - 1) x^{\alpha_1} y^{\alpha_2} z^{\alpha_3+1}. \end{aligned} \tag{38}$$

Pour tout $X \in A$, notons $D_0(1 \mapsto X)$ l'image par $D_0$ de l'application linéaire $\varphi_X \in K^1 = \mathrm{Hom}(V, A)$ définie par $\varphi_X(1) = X$. L'image de $D_0$ est constituée d'applications du type $D_0(1 \mapsto X)$, avec $X \in A$, et on a :

$$\begin{aligned} D_0(1 \mapsto X)(x) &= \; xX - Xx, \\ D_0(1 \mapsto X)(y) &= \; yX - Xy, \\ D_0(1 \mapsto X)(z) &= \; zX - Xz. \end{aligned}$$

On va réduire, modulo l'image de $D_0$, l'expression de $\varphi$ à :

$$\varphi(x) = \varphi(y) = 0, \; \varphi(z) = \sum_{\lambda^{\alpha_3 - 1} = 1} c_{\alpha_3} z^{\alpha_3}.$$

On concluera alors à l'aide du lemme 8.1.1 que les applications linéaires $\varphi$ définies par :

$$\varphi(x) = \varphi(y) = 0, \; \varphi(z) = z^{\alpha_3},$$

pour $\lambda^{\alpha_3 - 1} = 1$, forment une base de $HH^1(A)$, ce qui démontrera le lemme.
*Premier pas.* Identifions les termes de degré 0 en $z$ dans (37) :

$$\sum c_{\alpha_1, \alpha_2, 0}(1 - \lambda^1) x^{\alpha_1 + 1} y^{\alpha_2} + \sum c_{\alpha_1, \alpha_2, 0} \lambda^1 \alpha_2 x^{\alpha_1} y^{\alpha_2 - 1} = 0.$$

Pour tout $\alpha_2 \geq 0$, pour tout $\alpha_1 \geq 1$, on a :

$$c_{\alpha_1 - 1, \alpha_2, 0}(1 - \lambda) + c_{\alpha_1, \alpha_2 + 1, 0}(\alpha_2 + 1)\lambda = 0.$$

En raisonnant par récurrence descendante sur $\alpha_1$ pour $n = \alpha_1 - \alpha_2 \in \mathbb{Z}$ fixé, et puisque $\lambda \neq 1$ et que $\varphi(z)$ est un polynôme, on montre que tous les termes $c_{*,*,0}$ sont nuls. Ainsi on a montré que

$$\varphi(z) = \sum_{\alpha_3 \geq 1} c_\alpha x^{\alpha_1} y^{\alpha_2} z^{\alpha_3}.$$

*Deuxième pas.* Fixons $\alpha_3 \geq 1$, et $n = \alpha_2 - \alpha_1 \in \mathbb{Z}$ tel que $\lambda^n = 1$. Identifions alors le terme en $x^{\alpha_1 + 1} y^{\alpha_2} z^{\alpha_3}$ dans (37) :

$$c_\alpha(1 - \lambda^{1 - \alpha_3}) + c_{\alpha_1 + 1, \alpha_2 + 1, \alpha_3} \lambda^{1 - \alpha_3}(\alpha_2 + 1) = a_{\alpha_1 + 1, \alpha_2, \alpha_3 - 1}(\lambda^{1 - \alpha_1 - 1 + \alpha_2} - 1).$$

Puisqu'on est dans le cas où $\lambda^{\alpha_2 - \alpha_1} = 1$, le second terme de cette égalité est nul. Deux cas se présentent alors :

- si $\lambda^{1 - \alpha_3} \neq 1$, on obtient comme précédemment, par récurrence descendante, que $c_\alpha = 0$ pour tous $\alpha_1$ et $\alpha_2$ tels que $\alpha_2 - \alpha_1 = n$ ;
- si $\lambda^{1 - \alpha_3} = 1$, on obtient directement que $c_\alpha = 0$ pour tous $\alpha_1, \alpha_2 \geq 1$ tels que $\alpha_2 - \alpha_1 = n$.

*Troisième pas.* Posons alors

$$X_1 = \sum_{\substack{\alpha_1,\alpha_2,\alpha_3, \\ \alpha_3 \geq 1,\ \lambda^{\alpha_1-\alpha_2} \neq 1}} (\lambda^{-\alpha_2+\alpha_1} - 1)^{-1} c_\alpha x^{\alpha_1} y^{\alpha_2} z^{\alpha_3 - 1},$$

et remplaçons $\varphi$ par $\varphi - D_0(1 \mapsto X_1)$, ce qui ne change pas sa classe de cohomologie. Alors

$$\varphi(z) = \sum_{\alpha_3 \geq 1, \lambda^{\alpha_1-\alpha_2}=1} c_\alpha x^{\alpha_1} y^{\alpha_2} z^{\alpha_3},$$

c'est-à-dire, en tenant compte du pas précédent, que

$$\varphi(z) = \sum_{\lambda^{\alpha_3-1}=1} c_{\alpha_3} z^{\alpha_3}.$$

*Quatrième pas.* Le terme de gauche de (37) est alors identiquement nul. On en déduit que pour tous $\alpha_1, \alpha_2$ tels que $\lambda^{1-\alpha_1+\alpha_2} \neq 1$, alors $a_\alpha = 0$. Ainsi

$$\varphi(x) = \sum_{\lambda^{1-\alpha_1+\alpha_2}=1} a_\alpha x^{\alpha_1} y^{\alpha_2} z^{\alpha_3}.$$

Posons

$$X_2 = \sum_{\substack{\alpha_1,\alpha_2,\alpha_3, \\ \lambda^{1-\alpha_1+\alpha_2} = 1,\ \lambda^{-\alpha_3} = 1}} (\alpha_2 + 1)^{-1} a_\alpha x^{\alpha_1} y^{\alpha_2+1} z^{\alpha_3},$$

et remplaçons $\varphi$ par $\varphi - D_0(1 \mapsto X_2)$, ce qui ne change ni sa classe de cohomologie, ni sa valeur en $z$ par le lemme 8.1.1. On vérifie, toujours à l'aide du lemme 8.1.1, qu'on a alors

$$\varphi(x) = \sum_{\lambda^{1-\alpha_1+\alpha_2}=1,\ \lambda^{-\alpha_3}\neq 1} a_\alpha x^{\alpha_1} y^{\alpha_2} z^{\alpha_3}.$$

*Cinquième pas.* On se fixe $\alpha_3$ tel que $\lambda^{-\alpha_3} \neq 1$, et $n \in \mathbb{Z}$ tel que $\lambda^n = 1$. Notons $X(n)$ la partie de $\varphi(x)$ constituée des monômes $x^{\alpha_1} y^{\alpha_2} z^{\alpha_3}$ tels que $\alpha_2 = \alpha_1 + n - 1$. C'est un polynôme de la forme :

$$X(n) = \sum_{\alpha_1=0}^{p_n} a_\alpha x^{\alpha_1} y^{\alpha_1+n-1} z^{\alpha_3}.$$

Si $p_n \geq 1$, posons :

$$X_3^{n,p_n} = (1 - \lambda^{-\alpha_3})^{-1} a_{p_n, p_n+n-1, \alpha_3} x^{p_n - 1} y^{p_n+n-1} z^{\alpha_3}.$$

Remplaçons $\varphi$ par $\varphi - D_0(1 \mapsto X_3^{n,p_n})$, ce qui ne change ni sa classe de cohomologie, ni sa valeur en $z$ par le lemme 8.1.1. Par ailleurs on n'a changé dans $\varphi(x)$ que la partie $X(n)$, dont on a fait strictement baisser le degré. En raisonnant ainsi par récurrence descendante, et ceci pour tout $n$, on se ramène à :

$$\varphi(x) = \sum_{\lambda^{1+\alpha_2}=1, \lambda^{-\alpha_3}\neq 1} a_\alpha y^{\alpha_2} z^{\alpha_3}.$$

*Sixième pas.* Il découle de l'équation (38), dont le terme de gauche est nul, que dans $\varphi(y)$ n'apparaissent que des termes en $x^{\alpha_1}y^{\alpha_2}z^{\alpha_3}$ avec $\lambda^{-1-\alpha_1+\alpha_2} = 1$.
Etudions maintenant l'équation (36). Fixons $\alpha_3$ tel que $\lambda^{\alpha_3} = 1$. Identifions le terme en $x^{\alpha_1}y^{\alpha_2}z^{\alpha_3}$ :

$$b_{\alpha_1,\alpha_2+1,\alpha_3}(\alpha_2 + 1) = 0.$$

Ainsi $b_{\alpha_1,\alpha_2,\alpha_3} = 0$ dès que $\alpha_2 \geq 1$, et la partie ayant un exposant partiel en $z$ valant $\alpha_3$ dans $\varphi(y)$ vaut donc

$$Y(\alpha_3) = \sum_{\lambda^{-1-\alpha_1}=1} b_{\alpha_1,0,\alpha_3} x^{\alpha_1} z^{\alpha_3}.$$

Posons alors

$$X_4^{(\alpha_3)} = \sum_{\lambda^{-1-\alpha_1}=1} b_{\alpha_1,0,\alpha_3}(\alpha_1 + 1)^{-1} x^{\alpha_1+1} z^{\alpha_3},$$

et remplaçons $\varphi$ par $\varphi - D_0(1 \mapsto X_4^{(\alpha_3)})$, ce qui ne change ni sa classe de cohomologie, ni sa valeur en $z$, ni sa valeur en $x$ par le lemme 8.1.1. Ce même lemme montre également que $[y, X_4^{(\alpha_3)}] = Y(\alpha_3)$, donc en réitérant ceci pour tous les $\alpha_3$ tels que $\lambda^{\alpha_3} = 1$, on se ramène à :

$$\varphi(y) = \sum_{\lambda^{-1-\alpha_1+\alpha_2}=1, \lambda^{-\alpha_3}\neq 1} b_\alpha x^{\alpha_1} y^{\alpha_2} z^{\alpha_3}.$$

*Septième pas.* Soit $\alpha_3$ tel que $\lambda^{-\alpha_3} \neq 1$. Fixons $n \in \mathbb{Z}$ tel que $\lambda^n = 1$. Soient $\alpha_1 \geq 0, \alpha_2 \geq 1$ tels que $\alpha_2 = n + \alpha_1 + 1$, et identifions les termes en $x^{\alpha_1+1}y^{\alpha_2}z^{\alpha_3}$ dans (36) :

$$b_{\alpha_1,\alpha_2,\alpha_3}(1 - \lambda^{-\alpha_3}) + b_{\alpha_1+1,\alpha_2+1,\alpha_3}\lambda^{-\alpha_3}(\alpha_2 + 1) =$$
$$a_{\alpha_1+1,\alpha_2-1,\alpha_3}(1 - \lambda^{-\alpha_3}) - a_{\alpha_1+2,\alpha_2,\alpha_3}(\alpha_1 + 2).$$

Vue la forme à laquelle on s'est ramené pour $\varphi(x)$ au cinquième pas, le terme de droite de cette égalité est toujours nul. A nouveau par une récurrence descendante à $n$ fixé on en déduit que pour tout $\alpha_2 \geq 1$ on a : $b_{\alpha_1,\alpha_2,\alpha_3} = 0$. Ainsi

$$\varphi(y) = \sum_{\lambda^{-1-\alpha_1}=1, \lambda^{-\alpha_3}\neq 1} b_\alpha x^{\alpha_1} z^{\alpha_3}.$$

*Huitième pas.* La relation (36) devient alors :

$$\sum_{\lambda^{-1-\alpha_1}=1, \lambda^{-\alpha_3}\neq 1} b_\alpha x^{\alpha_1+1} z^{\alpha_3} = \sum_{\lambda^{1+\alpha_2}=1, \lambda^{-\alpha_3}\neq 1} a_\alpha y^{\alpha_2+1} z^{\alpha_3},$$

tous les $a_\alpha$ et les $b_\alpha$ sont donc nuls, et finalement $\varphi(x) = \varphi(y) = 0$.
On s'est ainsi ramené à la forme annoncée. Ainsi la famille d'applications linéaires $(\varphi_k)_{k\in\mathbb{N},\ \lambda^{1-k}=1}$ définies par :

$$\varphi_k(x) = \varphi_k(y) = 0,\ \varphi_k(z) = z^k,$$

engendre $\operatorname{Ker} D_1$ modulo $\operatorname{Im} D_0$. On vérifie que l'image de ces monômes est libre modulo $\operatorname{Im} D_0$ grâce au lemme 8.1.1 et à la liberté des $z^k$ dans $A$. □

### 8.2 Application à la question de la dualité pour l'homologie de $A$.

On peut maintenant énoncer le résultat suivant.

**Proposition 8.2.1** *Soit $\lambda \in k^*$, non racine de l'unité. Soit $A$ l'algèbre engendrée sur $k$ par $x, y, z$ avec les relations*

$$xy = yx + 1, \ xz = \lambda zx, \ yz = \lambda^{-1} zy.$$

*Alors il n'existe pas d'entier $d$ tel que pour tout $*$ on ait :*

$$HH_*(A) = HH^{d-*}(A).$$

*Preuve.* On a montré aux lemmes 8.1.2 et 8.1.3 que $HH^0(A) = HH^1(A) = k$. Si $A$ vérifiait une dualité du type $HH_*(A) = HH^{d-*}(A)$ pour un certain entier $d$, on devrait donc avoir $HH_d(A) = HH_{d-1}(A) = k$, ce qui est impossible d'après le théorème 5.5.2. □

## Remerciements.

Nous tenons à remercier Cyrille Ospel et Thierry Lambre pour des discussions fructueuses au cours de l'élaboration de ce travail.

## Références